\documentstyle[amsfonts,thmsk,11pt,a4,sw20lart]{article}
\input tcilatex

\makeatletter

\def\diagram{\leftwidth=\z@ \rightwidth=\z@ \topheight=\z@
\botheight=\z@ \setbox\@picbox\hbox\bgroup}

\def\enddiagram{\egroup\wd\@picbox\rightwidth\unitlength
\ht\@picbox\topheight\unitlength \dp\@picbox\botheight\unitlength
\hskip\leftwidth\unitlength\box\@picbox}

\def\bfig{\begin{diagram}}
\def\efig{\end{diagram}}
\newcount\wideness \newcount\leftwidth \newcount\rightwidth
\newcount\highness \newcount\topheight \newcount\botheight

\def\ratchet#1#2{\ifnum#1<#2 \global #1=#2 \fi}

\def\putbox(#1,#2)#3{%
\horsize{\wideness}{#3} \divide\wideness by 2
{\advance\wideness by #1 \ratchet{\rightwidth}{\wideness}}
{\advance\wideness by -#1 \ratchet{\leftwidth}{\wideness}}
\vertsize{\highness}{#3} \divide\highness by 2
{\advance\highness by #2 \ratchet{\topheight}{\highness}}
{\advance\highness by -#2 \ratchet{\botheight}{\highness}}
\put(#1,#2){\makebox(0,0){$#3$}}}

\def\putlbox(#1,#2)#3{%
\horsize{\wideness}{#3}
{\advance\wideness by #1 \ratchet{\rightwidth}{\wideness}}
{\ratchet{\leftwidth}{-#1}}
\vertsize{\highness}{#3} \divide\highness by 2
{\advance\highness by #2 \ratchet{\topheight}{\highness}}
{\advance\highness by -#2 \ratchet{\botheight}{\highness}}
\put(#1,#2){\makebox(0,0)[l]{$#3$}}}

\def\putrbox(#1,#2)#3{%
\horsize{\wideness}{#3}
{\ratchet{\rightwidth}{#1}}
{\advance\wideness by -#1 \ratchet{\leftwidth}{\wideness}}
\vertsize{\highness}{#3} \divide\highness by 2
{\advance\highness by #2 \ratchet{\topheight}{\highness}}
{\advance\highness by -#2 \ratchet{\botheight}{\highness}}
\put(#1,#2){\makebox(0,0)[r]{$#3$}}}

\def\adjust[#1]{} % For compatibility

\newcount \coefa
\newcount \coefb
\newcount \coefc
\newcount\tempcounta
\newcount\tempcountb
\newcount\tempcountc
\newcount\tempcountd
\newcount\xext
\newcount\yext
\newcount\xoff
\newcount\yoff
\newcount\gap%
\newcount\arrowtypea
\newcount\arrowtypeb
\newcount\arrowtypec
\newcount\arrowtyped
\newcount\arrowtypee
\newcount\height
\newcount\width
\newcount\xpos
\newcount\ypos
\newcount\run
\newcount\rise
\newcount\arrowlength
\newcount\halflength
\newcount\arrowtype
\newdimen\tempdimen
\newdimen\xlen
\newdimen\ylen
\newsavebox{\tempboxa}%
\newsavebox{\tempboxb}%
\newsavebox{\tempboxc}%

\newdimen\w@dth

\def\setw@dth#1#2{\setbox\z@\hbox{$#1$}\w@dth=\wd\z@
\setbox\@ne\hbox{$#2$}\ifnum\w@dth<\wd\@ne \w@dth=\wd\@ne \fi
\advance\w@dth by 1.2em}

\def\t@^#1_#2{\def\n@one{#1}\def\n@two{#2}\mathrel{\setw@dth{#1}{#2}
\mathop{\hbox to \w@dth{\rightarrowfill}}\limits
\ifx\n@one\empty\else ^{\box\z@}\fi
\ifx\n@two\empty\else _{\box\@ne}\fi}}
\def\t@@^#1{\@ifnextchar_ {\t@^{#1}}{\t@^{#1}_{}}}
\def\to{\@ifnextchar^ {\t@@}{\t@@^{}}}

\def\t@left^#1_#2{\def\n@one{#1}\def\n@two{#2}\mathrel{\setw@dth{#1}{#2}
\mathop{\hbox to \w@dth{\leftarrowfill}}\limits
\ifx\n@one\empty\else ^{\box\z@}\fi
\ifx\n@two\empty\else _{\box\@ne}\fi}}
\def\t@@left^#1{\@ifnextchar_ {\t@left^{#1}}{\t@left^{#1}_{}}}
\def\toleft{\@ifnextchar^ {\t@@left}{\t@@left^{}}}

\def\two@^#1_#2{\def\n@one{#1}\def\n@two{#2}\mathrel{\setw@dth{#1}{#2}
\mathop{\vcenter{\hbox to \w@dth{\rightarrowfill}\kern-1.7ex
                 \hbox to \w@dth{\rightarrowfill}}%
       }\limits
\ifx\n@one\empty\else ^{\box\z@}\fi
\ifx\n@two\empty\else _{\box\@ne}\fi}}
\def\tw@@^#1{\@ifnextchar_ {\two@^{#1}}{\two@^{#1}_{}}}
\def\two{\@ifnextchar^ {\tw@@}{\tw@@^{}}}

\def\tofr@^#1_#2{\def\n@one{#1}\def\n@two{#2}\mathrel{\setw@dth{#1}{#2}
\mathop{\vcenter{\hbox to \w@dth{\rightarrowfill}\kern-1.7ex
                 \hbox to \w@dth{\leftarrowfill}}%
       }\limits
\ifx\n@one\empty\else ^{\box\z@}\fi
\ifx\n@two\empty\else _{\box\@ne}\fi}}
\def\t@fr@^#1{\@ifnextchar_ {\tofr@^{#1}}{\tofr@^{#1}_{}}}
\def\tofro{\@ifnextchar^ {\t@fr@}{\t@fr@^{}}}

\def\mon{\mathop{\m@th\hbox to
      14.6\P@{\lasyb\char'51\hskip-2.1\P@$\arrext$\hss
$\mathord\rightarrow$}}\limits} % width of \epi
\def\leftmono{\mathrel{\m@th\hbox to
14.6\P@{$\mathord\leftarrow$\hss$\arrext$\hskip-2.1\P@\lasyb\char'50%
}}\limits} % width of \epi
\mathchardef\arrext="0200       % amr minus for arrow extension (see \into)

\def\settypes(#1,#2,#3){\arrowtypea#1 \arrowtypeb#2 \arrowtypec#3}
\def\settoheight#1#2{\setbox\@tempboxa\hbox{#2}#1\ht\@tempboxa\relax}%
\def\settodepth#1#2{\setbox\@tempboxa\hbox{#2}#1\dp\@tempboxa\relax}%
\def\settokens[#1`#2`#3`#4]{%
     \def\tokena{#1}\def\tokenb{#2}\def\tokenc{#3}\def\tokend{#4}}
\def\setsqparms[#1`#2`#3`#4;#5`#6]{%
\arrowtypea #1
\arrowtypeb #2
\arrowtypec #3
\arrowtyped #4
\width #5
\height #6
}
\def\setpos(#1,#2){\xpos=#1 \ypos#2}

\def\settriparms[#1`#2`#3;#4]{\settripairparms[#1`#2`#3`1`1;#4]}%

\def\settripairparms[#1`#2`#3`#4`#5;#6]{%
\arrowtypea #1
\arrowtypeb #2
\arrowtypec #3
\arrowtyped #4
\arrowtypee #5
\width #6
\height #6
}

\def\resetparms{\settripairparms[1`1`1`1`1;500]\width 500}%default values%

\resetparms

\def\mvector(#1,#2)#3{%%
\put(0,0){\vector(#1,#2){#3}}%
\put(0,0){\vector(#1,#2){26}}%
}
\def\evector(#1,#2)#3{{%%
\arrowlength #3
\put(0,0){\vector(#1,#2){\arrowlength}}%
\advance \arrowlength by-30
\put(0,0){\vector(#1,#2){\arrowlength}}%
}}

\def\horsize#1#2{%
\settowidth{\tempdimen}{$#2$}%
#1=\tempdimen
\divide #1 by\unitlength
}

\def\vertsize#1#2{%
\settoheight{\tempdimen}{$#2$}%
#1=\tempdimen
\settodepth{\tempdimen}{$#2$}%
\advance #1 by\tempdimen
\divide #1 by\unitlength
}

\def\putvector(#1,#2)(#3,#4)#5#6{{%
\ifnum3<\arrowtype
\putdashvector(#1,#2)(#3,#4)#5\arrowtype
\else
\ifnum\arrowtype<-3
\putdashvector(#1,#2)(#3,#4)#5\arrowtype
\else
\xpos=#1
\ypos=#2
\run=#3
\rise=#4
\arrowlength=#5
\ifnum \arrowtype<0
    \ifnum \run=0
        \advance \ypos by-\arrowlength
    \else
        \tempcounta \arrowlength
        \multiply \tempcounta by\rise
        \divide \tempcounta by\run
        \ifnum\run>0
            \advance \xpos by\arrowlength
            \advance \ypos by\tempcounta
        \else
            \advance \xpos by-\arrowlength
            \advance \ypos by-\tempcounta
        \fi
    \fi
    \multiply \arrowtype by-1
    \multiply \rise by-1
    \multiply \run by-1
\fi
\ifcase \arrowtype
\or \put(\xpos,\ypos){\vector(\run,\rise){\arrowlength}}%
\or \put(\xpos,\ypos){\mvector(\run,\rise)\arrowlength}%
\or \put(\xpos,\ypos){\evector(\run,\rise){\arrowlength}}%
\fi\fi\fi
}}

\def\putsplitvector(#1,#2)#3#4{%%
\xpos #1
\ypos #2
\arrowtype #4
\halflength #3
\arrowlength #3
\gap 140
\advance \halflength by-\gap
\divide \halflength by2
\ifnum\arrowtype>0
   \ifcase \arrowtype
   \or \put(\xpos,\ypos){\line(0,-1){\halflength}}%
       \advance\ypos by-\halflength
       \advance\ypos by-\gap
       \put(\xpos,\ypos){\vector(0,-1){\halflength}}%
   \or \put(\xpos,\ypos){\line(0,-1)\halflength}%
       \put(\xpos,\ypos){\vector(0,-1)3}%
       \advance\ypos by-\halflength
       \advance\ypos by-\gap
       \put(\xpos,\ypos){\vector(0,-1){\halflength}}%
   \or \put(\xpos,\ypos){\line(0,-1)\halflength}%
       \advance\ypos by-\halflength
       \advance\ypos by-\gap
       \put(\xpos,\ypos){\evector(0,-1){\halflength}}%
   \fi
\else \arrowtype=-\arrowtype
   \ifcase\arrowtype
   \or \advance \ypos by-\arrowlength
       \put(\xpos,\ypos){\line(0,1){\halflength}}%
       \advance\ypos by\halflength
       \advance\ypos by\gap
       \put(\xpos,\ypos){\vector(0,1){\halflength}}%
   \or \advance \ypos by-\arrowlength
       \put(\xpos,\ypos){\line(0,1)\halflength}%
       \put(\xpos,\ypos){\vector(0,1)3}%
       \advance\ypos by\halflength
       \advance\ypos by\gap
       \put(\xpos,\ypos){\vector(0,1){\halflength}}%
   \or \advance \ypos by-\arrowlength
       \put(\xpos,\ypos){\line(0,1)\halflength}%
       \advance\ypos by\halflength
       \advance\ypos by\gap
       \put(\xpos,\ypos){\evector(0,1){\halflength}}%
   \fi
\fi
}

\def\putmorphism(#1)(#2,#3)[#4`#5`#6]#7#8#9{{%
\run #2
\rise #3
\ifnum\rise=0
  \puthmorphism(#1)[#4`#5`#6]{#7}{#8}#9%
\else\ifnum\run=0
  \putvmorphism(#1)[#4`#5`#6]{#7}{#8}#9%
\else
\setpos(#1)%
\arrowlength #7
\arrowtype #8
\ifnum\run=0
\else\ifnum\rise=0
\else
\ifnum\run>0
    \coefa=1
\else
   \coefa=-1
\fi
\ifnum\arrowtype>0
   \coefb=0
   \coefc=-1
\else
   \coefb=\coefa
   \coefc=1
   \arrowtype=-\arrowtype
\fi
\width=2
\multiply \width by\run
\divide \width by\rise
\ifnum \width<0  \width=-\width\fi
\advance\width by60
\if l#9 \width=-\width\fi
\putbox(\xpos,\ypos){#4}%            %node 1
{\multiply \coefa by\arrowlength%      %node 2
\advance\xpos by\coefa
\multiply \coefa by\rise
\divide \coefa by\run
\advance \ypos by\coefa
\putbox(\xpos,\ypos){#5} }%
{\multiply \coefa by\arrowlength%      %label
\divide \coefa by2
\advance \xpos by\coefa
\advance \xpos by\width
\multiply \coefa by\rise
\divide \coefa by\run
\advance \ypos by\coefa
\if l#9%
   \putrbox(\xpos,\ypos){#6}%
\else\if r#9%
   \putlbox(\xpos,\ypos){#6}%
\fi\fi }%
{\multiply \rise by-\coefc%             %arrow
\multiply \run by-\coefc
\multiply \coefb by\arrowlength
\advance \xpos by\coefb
\multiply \coefb by\rise
\divide \coefb by\run
\advance \ypos by\coefb
\multiply \coefc by70
\advance \ypos by\coefc
\multiply \coefc by\run
\divide \coefc by\rise
\advance \xpos by\coefc
\multiply \coefa by140
\multiply \coefa by\run
\divide \coefa by\rise
\advance \arrowlength by\coefa
\ifcase\arrowtype
\or \put(\xpos,\ypos){\vector(\run,\rise){\arrowlength}}%
\or \put(\xpos,\ypos){\mvector(\run,\rise){\arrowlength}}%
\or \put(\xpos,\ypos){\evector(\run,\rise){\arrowlength}}%
\fi}\fi\fi\fi\fi}}

\newcount\numbdashes \newcount\lengthdash \newcount\increment

\def\howmanydashes{% Actually returns both number and length
\numbdashes=\arrowlength \lengthdash=40
\divide\numbdashes by \lengthdash
\lengthdash=\arrowlength
\divide\lengthdash by \numbdashes
\increment=\lengthdash
\multiply\lengthdash by 3
\divide\lengthdash by 5
}

\def\putdashvector(#1)(#2,#3)#4#5{%
\ifnum#3=0 \putdashhvector(#1){#4}#5
\else
\ifnum#2=0
\putdashvvector(#1){#4}#5\fi\fi}

\def\putdashhvector(#1,#2)#3#4{{%
\arrowlength=#3 \howmanydashes
\multiput(#1,#2)(\increment,0){\numbdashes}%
{\vrule height .4pt width \lengthdash\unitlength}
\arrowtype=#4 \xpos=#1
\ifnum\arrowtype<0 \advance\arrowtype by 7 \fi
\ifcase\arrowtype
\or \advance\xpos by 10
    \put(\xpos,#2){\vector(-1,0){\lengthdash}}
    \advance\xpos by 40
    \put(\xpos,#2){\vector(-1,0){\lengthdash}}
\or \advance \xpos by 10
    \put(\xpos,#2){\vector(-1,0){\lengthdash}}
    \advance\xpos by  \arrowlength
    \advance\xpos by  -50
    \put(\xpos,#2){\vector(-1,0){\lengthdash}}
\or \advance\xpos by 10
    \put(\xpos,#2){\vector(-1,0){\lengthdash}}
\or \advance\xpos by \arrowlength
    \advance\xpos by -\lengthdash
    \put(\xpos,#2){\vector(1,0){\lengthdash}}
\or {\advance\xpos by 10
    \put(\xpos,#2){\vector(1,0){\lengthdash}}}
    \advance\xpos by \arrowlength
    \advance\xpos by -\lengthdash
    \put(\xpos,#2){\vector(1,0){\lengthdash}}
\or \advance\xpos by \arrowlength
    \advance\xpos by -\lengthdash
    \put(\xpos,#2){\vector(1,0){\lengthdash}}
    \advance\xpos by -40
    \put(\xpos,#2){\vector(1,0){\lengthdash}}
   \fi
}}

\def\putdashvvector(#1,#2)#3#4{{%
\arrowlength=#3 \howmanydashes
\ypos=#2 \advance\ypos by -\arrowlength
\multiput(#1,#2)(0,\increment){\numbdashes}%
    {\vrule width .4pt height \lengthdash\unitlength}
\arrowtype=#4 \ypos=#2
\ifnum\arrowtype<0 \advance\arrowtype by 7 \fi
\ifcase\arrowtype
\or \advance\ypos by \arrowlength \advance\ypos by -40
    \put(#1,\ypos){\vector(0,1){\lengthdash}}
    \advance\ypos by -40
    \put(#1,\ypos){\vector(0,1){\lengthdash}}
\or \advance\ypos by 10
    \put(#1,\ypos){\vector(0,1){\lengthdash}}
    \advance\ypos by \arrowlength \advance\ypos by -40
    \put(#1,\ypos){\vector(0,1){\lengthdash}}
\or \advance\ypos by \arrowlength \advance\ypos by -40
    \put(#1,\ypos){\vector(0,1){\lengthdash}}
\or \advance\ypos by 10
    \put(#1,\ypos){\vector(0,-1){\lengthdash}}
\or \advance\ypos by 10
    \put(#1,\ypos){\vector(0,-1){\lengthdash}}
    \advance\ypos by \arrowlength \advance\ypos by -40
    \put(#1,\ypos){\vector(0,-1){\lengthdash}}
\or \advance\ypos by 10
    \put(#1,\ypos){\vector(0,-1){\lengthdash}}
    \advance\ypos by 40
    \put(#1,\ypos){\vector(0,-1){\lengthdash}}
\fi
}}

\def\puthmorphism(#1,#2)[#3`#4`#5]#6#7#8{{%
\xpos #1
\ypos #2
\width #6
\arrowlength #6
\arrowtype=#7
\putbox(\xpos,\ypos){#3\vphantom{#4}}%
{\advance \xpos by\arrowlength
\putbox(\xpos,\ypos){\vphantom{#3}#4}}%
\horsize{\tempcounta}{#3}%
\horsize{\tempcountb}{#4}%
\divide \tempcounta by2
\divide \tempcountb by2
\advance \tempcounta by30
\advance \tempcountb by30
\advance \xpos by\tempcounta
\advance \arrowlength by-\tempcounta
\advance \arrowlength by-\tempcountb
\putvector(\xpos,\ypos)(1,0)\arrowlength\arrowtype
\divide \arrowlength by2
\advance \xpos by\arrowlength
\vertsize{\tempcounta}{#5}%
\divide\tempcounta by2
\advance \tempcounta by20
\if a#8 %
   \advance \ypos by\tempcounta
   \putbox(\xpos,\ypos){#5}%
\else
   \advance \ypos by-\tempcounta
   \putbox(\xpos,\ypos){#5}%
\fi}}

\def\putvmorphism(#1,#2)[#3`#4`#5]#6#7#8{{%
\xpos #1
\ypos #2
\arrowlength #6
\arrowtype #7
\settowidth{\xlen}{$#5$}%
\putbox(\xpos,\ypos){#3}%
{\advance \ypos by-\arrowlength
\putbox(\xpos,\ypos){#4}}%
{\advance\arrowlength by-140
\advance \ypos by-70
\ifdim\xlen>0pt
   \if m#8%
      \putsplitvector(\xpos,\ypos)\arrowlength\arrowtype
   \else
   \putvector(\xpos,\ypos)(0,-1)\arrowlength\arrowtype
   \fi
\else
   \putvector(\xpos,\ypos)(0,-1)\arrowlength\arrowtype
\fi}%
\ifdim\xlen>0pt
   \divide \arrowlength by2
   \advance\ypos by-\arrowlength
   \if l#8%
      \advance \xpos by-40
      \putrbox(\xpos,\ypos){#5}%
   \else\if r#8%
      \advance \xpos by40
      \putlbox(\xpos,\ypos){#5}%
   \else
      \putbox(\xpos,\ypos){#5}%
   \fi\fi
\fi
}}

\def\putsquarep<#1>(#2)[#3;#4`#5`#6`#7]{{%
\setsqparms[#1]%
\setpos(#2)%
\settokens[#3]%
\puthmorphism(\xpos,\ypos)[\tokenc`\tokend`{#7}]{\width}{\arrowtyped}b%
\advance\ypos by \height
\puthmorphism(\xpos,\ypos)[\tokena`\tokenb`{#4}]{\width}{\arrowtypea}a%
\putvmorphism(\xpos,\ypos)[``{#5}]{\height}{\arrowtypeb}l%
\advance\xpos by \width
\putvmorphism(\xpos,\ypos)[``{#6}]{\height}{\arrowtypec}r%
}}

\def\putsquare{\@ifnextchar <{\putsquarep}{\putsquarep%
   <\arrowtypea`\arrowtypeb`\arrowtypec`\arrowtyped;\width`\height>}}
\def\square{\@ifnextchar< {\squarep}{\squarep
   <\arrowtypea`\arrowtypeb`\arrowtypec`\arrowtyped;\width`\height>}}
\def\squarep<#1>[#2`#3`#4`#5;#6`#7`#8`#9]{{%       %     #2------>#3
\setsqparms[#1]%                                   %      |       |
\diagram%                                          %      |       |
\putsquarep<\arrowtypea`\arrowtypeb`\arrowtypec`%  %    #7|       |#8
\arrowtyped;\width`\height>%                       %      |       |
(0,0)[#2`#3`#4`{#5};#6`#7`#8`{#9}]%                %      |       |
\enddiagram%                                       %      v       v
}}                                                 %     #4------>#5
\def\putptrianglep<#1>(#2,#3)[#4`#5`#6;#7`#8`#9]{{%
\settriparms[#1]%
\xpos=#2 \ypos=#3
\advance\ypos by \height
\puthmorphism(\xpos,\ypos)[#4`#5`{#7}]{\height}{\arrowtypea}a%
\putvmorphism(\xpos,\ypos)[`#6`{#8}]{\height}{\arrowtypeb}l%
\advance\xpos by\height
\putmorphism(\xpos,\ypos)(-1,-1)[``{#9}]{\height}{\arrowtypec}r%
}}

\def\putptriangle{\@ifnextchar <{\putptrianglep}{\putptrianglep
   <\arrowtypea`\arrowtypeb`\arrowtypec;\height>}}
\def\ptriangle{\@ifnextchar <{\ptrianglep}{\ptrianglep
   <\arrowtypea`\arrowtypeb`\arrowtypec;\height>}}
\def\ptrianglep<#1>[#2`#3`#4;#5`#6`#7]{{%%    %      #2----->#3
\settriparms[#1]%                             %      |      /
\diagram%                                     %      |     /
\putptrianglep<\arrowtypea`\arrowtypeb`%      %    #6|    /#7
\arrowtypec;\height>%                         %      |   /
(0,0)[#2`#3`#4;#5`#6`{#7}]%                   %      |  /
\enddiagram%%                                 %      v v
}}                                            %      #4

\def\putqtrianglep<#1>(#2,#3)[#4`#5`#6;#7`#8`#9]{{%
\settriparms[#1]%
\xpos=#2 \ypos=#3
\advance\ypos by\height
\puthmorphism(\xpos,\ypos)[#4`#5`{#7}]{\height}{\arrowtypea}a%
\putmorphism(\xpos,\ypos)(1,-1)[``{#8}]{\height}{\arrowtypeb}l%
\advance\xpos by\height
\putvmorphism(\xpos,\ypos)[`#6`{#9}]{\height}{\arrowtypec}r%
}}

\def\putqtriangle{\@ifnextchar <{\putqtrianglep}{\putqtrianglep
   <\arrowtypea`\arrowtypeb`\arrowtypec;\height>}}
\def\qtriangle{\@ifnextchar <{\qtrianglep}{\qtrianglep
   <\arrowtypea`\arrowtypeb`\arrowtypec;\height>}}
\def\qtrianglep<#1>[#2`#3`#4;#5`#6`#7]{{%%    %        #2----->#3
\settriparms[#1]%                             %         \      |
\width=\height                                %          \     |
\diagram%                                     %         #6\    |#7
\putqtrianglep<\arrowtypea`\arrowtypeb`%      %            \   |
\arrowtypec;\height>%                         %             \  |
(0,0)[#2`#3`#4;#5`#6`{#7}]%                   %              v v
\enddiagram%%                                 %               #4
}}

\def\putdtrianglep<#1>(#2,#3)[#4`#5`#6;#7`#8`#9]{{%
\settriparms[#1]%
\xpos=#2 \ypos=#3
\puthmorphism(\xpos,\ypos)[#5`#6`{#9}]{\height}{\arrowtypec}b%
\advance\xpos by \height \advance\ypos by\height
\putmorphism(\xpos,\ypos)(-1,-1)[``{#7}]{\height}{\arrowtypea}l%
\putvmorphism(\xpos,\ypos)[#4``{#8}]{\height}{\arrowtypeb}r%
}}

\def\putdtriangle{\@ifnextchar <{\putdtrianglep}{\putdtrianglep
   <\arrowtypea`\arrowtypeb`\arrowtypec;\height>}}
\def\dtriangle{\@ifnextchar <{\dtrianglep}{\dtrianglep
   <\arrowtypea`\arrowtypeb`\arrowtypec;\height>}}
\def\dtrianglep<#1>[#2`#3`#4;#5`#6`#7]{{%%    %                  / |
\settriparms[#1]%                             %                 /  |
\width=\height                                %              #5/   |#6
\diagram%                                     %               /    |
\putdtrianglep<\arrowtypea`\arrowtypeb`%      %              /     |
\arrowtypec;\height>%                         %             v      v
(0,0)[#2`#3`#4;#5`#6`{#7}]%                   %            #3----->#4
\enddiagram%%                                 %                #7
}}

\def\putbtrianglep<#1>(#2,#3)[#4`#5`#6;#7`#8`#9]{{%
\settriparms[#1]%
\xpos=#2 \ypos=#3
\puthmorphism(\xpos,\ypos)[#5`#6`{#9}]{\height}{\arrowtypec}b%
\advance\ypos by\height
\putmorphism(\xpos,\ypos)(1,-1)[``{#8}]{\height}{\arrowtypeb}r%
\putvmorphism(\xpos,\ypos)[#4``{#7}]{\height}{\arrowtypea}l%
}}

\def\putbtriangle{\@ifnextchar <{\putbtrianglep}{\putbtrianglep
   <\arrowtypea`\arrowtypeb`\arrowtypec;\height>}}
\def\btriangle{\@ifnextchar <{\btrianglep}{\btrianglep
   <\arrowtypea`\arrowtypeb`\arrowtypec;\height>}}
\def\btrianglep<#1>[#2`#3`#4;#5`#6`#7]{{%%   %              | \
\settriparms[#1]%                            %              |  \
\width=\height                               %            #5|   \#6
\diagram%                                    %              |    \
\putbtrianglep<\arrowtypea`\arrowtypeb`%     %              |     \
\arrowtypec;\height>%                        %              v      v
(0,0)[#2`#3`#4;#5`#6`{#7}]%                  %              #3----->#4
\enddiagram%%                                %                 #7
}}

\def\putAtrianglep<#1>(#2,#3)[#4`#5`#6;#7`#8`#9]{{%
\settriparms[#1]%
\xpos=#2 \ypos=#3
{\multiply \height by2
\puthmorphism(\xpos,\ypos)[#5`#6`{#9}]{\height}{\arrowtypec}b}%
\advance\xpos by\height \advance\ypos by\height
\putmorphism(\xpos,\ypos)(-1,-1)[#4``{#7}]{\height}{\arrowtypea}l%
\putmorphism(\xpos,\ypos)(1,-1)[``{#8}]{\height}{\arrowtypeb}r%
}}

\def\putAtriangle{\@ifnextchar <{\putAtrianglep}{\putAtrianglep
   <\arrowtypea`\arrowtypeb`\arrowtypec;\height>}}
\def\Atriangle{\@ifnextchar <{\Atrianglep}{\Atrianglep
   <\arrowtypea`\arrowtypeb`\arrowtypec;\height>}}
\def\Atrianglep<#1>[#2`#3`#4;#5`#6`#7]{{%%         %         /   \
\settriparms[#1]%                                  %        /     \
\width=\height                                     %     #5/       \#6
\diagram%                                          %      /         \
\putAtrianglep<\arrowtypea`\arrowtypeb`%           %     /           \
\arrowtypec;\height>%                              %    v             v
(0,0)[#2`#3`#4;#5`#6`{#7}]%                        %   #3------------>#4
\enddiagram%%                                      %          #7
}}

\def\putAtrianglepairp<#1>(#2)[#3;#4`#5`#6`#7`#8]{{%
\settripairparms[#1]%
\setpos(#2)%
\settokens[#3]%
\puthmorphism(\xpos,\ypos)[\tokenb`\tokenc`{#7}]{\height}{\arrowtyped}b%
\advance\xpos by\height
\puthmorphism(\xpos,\ypos)[\phantom{\tokenc}`\tokend`{#8}]%
{\height}{\arrowtypee}b%
\advance\ypos by\height
\putmorphism(\xpos,\ypos)(-1,-1)[\tokena``{#4}]{\height}{\arrowtypea}l%
\putvmorphism(\xpos,\ypos)[``{#5}]{\height}{\arrowtypeb}m%
\putmorphism(\xpos,\ypos)(1,-1)[``{#6}]{\height}{\arrowtypec}r%
}}

\def\putAtrianglepair{\@ifnextchar <{\putAtrianglepairp}{\putAtrianglepairp%
   <\arrowtypea`\arrowtypeb`\arrowtypec`\arrowtyped`\arrowtypee;\height>}}
\def\Atrianglepair{\@ifnextchar <{\Atrianglepairp}{\Atrianglepairp%
   <\arrowtypea`\arrowtypeb`\arrowtypec`\arrowtyped`\arrowtypee;\height>}}

\def\Atrianglepairp<#1>[#2;#3`#4`#5`#6`#7]{{%           %  #2a
\settripairparms[#1]%                         %           / | \
\settokens[#2]%                               %          /  |  \
\width=\height                                %       #3/  #4   \#5
\diagram%                                     %        /    |    \
\putAtrianglepairp                            %       /     |     \
<\arrowtypea`\arrowtypeb`\arrowtypec`%        %      v      v      v
\arrowtyped`\arrowtypee;\height>%             %     #2b---->#2c---->#2d
(0,0)[{#2};#3`#4`#5`#6`{#7}]%                 %         #6     #7
\enddiagram%%
}}

\def\putVtrianglep<#1>(#2,#3)[#4`#5`#6;#7`#8`#9]{{%
\settriparms[#1]%
\xpos=#2 \ypos=#3
\advance\ypos by\height
{\multiply\height by2
\puthmorphism(\xpos,\ypos)[#4`#5`{#7}]{\height}{\arrowtypea}a}%
\putmorphism(\xpos,\ypos)(1,-1)[`#6`{#8}]{\height}{\arrowtypeb}l%
\advance\xpos by\height
\advance\xpos by\height
\putmorphism(\xpos,\ypos)(-1,-1)[``{#9}]{\height}{\arrowtypec}r%
}}

\def\putVtriangle{\@ifnextchar <{\putVtrianglep}{\putVtrianglep
   <\arrowtypea`\arrowtypeb`\arrowtypec;\height>}}
\def\Vtriangle{\@ifnextchar <{\Vtrianglep}{\Vtrianglep
   <\arrowtypea`\arrowtypeb`\arrowtypec;\height>}}
\def\Vtrianglep<#1>[#2`#3`#4;#5`#6`#7]{{%%     %        #2------------->#3
\settriparms[#1]%                              %         \             /
\width=\height                                 %          \           /
\diagram%                                      %         #6\         /#7
\putVtrianglep<\arrowtypea`\arrowtypeb`%       %            \       /
\arrowtypec;\height>%                          %             \     /
(0,0)[#2`#3`#4;#5`#6`{#7}]%                    %              v   v
\enddiagram%%                                  %               #4
}}

\def\putVtrianglepairp<#1>(#2)[#3;#4`#5`#6`#7`#8]{{
\settripairparms[#1]%
\setpos(#2)%
\settokens[#3]%
\advance\ypos by\height
\putmorphism(\xpos,\ypos)(1,-1)[`\tokend`{#6}]{\height}{\arrowtypec}l%
\puthmorphism(\xpos,\ypos)[\tokena`\tokenb`{#4}]{\height}{\arrowtypea}a%
\advance\xpos by\height
\puthmorphism(\xpos,\ypos)[\phantom{\tokenb}`\tokenc`{#5}]%
{\height}{\arrowtypeb}a%
\putvmorphism(\xpos,\ypos)[``{#7}]{\height}{\arrowtyped}m%
\advance\xpos by\height
\putmorphism(\xpos,\ypos)(-1,-1)[``{#8}]{\height}{\arrowtypee}r%
}}

\def\putVtrianglepair{\@ifnextchar <{\putVtrianglepairp}{\putVtrianglepairp%
    <\arrowtypea`\arrowtypeb`\arrowtypec`\arrowtyped`\arrowtypee;\height>}}
\def\Vtrianglepair{\@ifnextchar <{\Vtrianglepairp}{\Vtrianglepairp%
    <\arrowtypea`\arrowtypeb`\arrowtypec`\arrowtyped`\arrowtypee;\height>}}
\def\Vtrianglepairp<#1>[#2;#3`#4`#5`#6`#7]{{%  %  #2a---->#2b---->#2c
\settripairparms[#1]%                          %   \      |      /
\settokens[#2]%                                %    \     |     /
\diagram%                                      %   #5\   #6    /#7
\putVtrianglepairp                             %      \   |   /
<\arrowtypea`\arrowtypeb`\arrowtypec`%         %       \  |  /
\arrowtyped`\arrowtypee;\height>%              %        v v v
(0,0)[{#2};#3`#4`#5`#6`{#7}]%                  %         #2d
\enddiagram%%
}}

\def\putCtrianglep<#1>(#2,#3)[#4`#5`#6;#7`#8`#9]{{%
\settriparms[#1]%
\xpos=#2 \ypos=#3
\advance\ypos by\height
\putmorphism(\xpos,\ypos)(1,-1)[``{#9}]{\height}{\arrowtypec}l%
\advance\xpos by\height
\advance\ypos by\height
\putmorphism(\xpos,\ypos)(-1,-1)[#4`#5`{#7}]{\height}{\arrowtypea}l%
{\multiply\height by 2
\putvmorphism(\xpos,\ypos)[`#6`{#8}]{\height}{\arrowtypeb}r}%
}}

\def\putCtriangle{\@ifnextchar <{\putCtrianglep}{\putCtrianglep
    <\arrowtypea`\arrowtypeb`\arrowtypec;\height>}}
\def\Ctriangle{\@ifnextchar <{\Ctrianglep}{\Ctrianglep
    <\arrowtypea`\arrowtypeb`\arrowtypec;\height>}}
\def\Ctrianglep<#1>[#2`#3`#4;#5`#6`#7]{{%%   %                / |
\settriparms[#1]%                            %             #5/  |
\width=\height                               %              /   |
\diagram%                                    %             v    |
\putCtrianglep<\arrowtypea`\arrowtypeb`%     %           #3     |#6
\arrowtypec;\height>%                        %             \    |
(0,0)[#2`#3`#4;#5`#6`{#7}]%                  %            #7\   |
\enddiagram%%                                %               \  |
}}                                           %                v v
\def\putDtrianglep<#1>(#2,#3)[#4`#5`#6;#7`#8`#9]{{%
\settriparms[#1]%
\xpos=#2 \ypos=#3
\advance\xpos by\height \advance\ypos by\height
\putmorphism(\xpos,\ypos)(-1,-1)[``{#9}]{\height}{\arrowtypec}r%
\advance\xpos by-\height \advance\ypos by\height
\putmorphism(\xpos,\ypos)(1,-1)[`#5`{#8}]{\height}{\arrowtypeb}r%
{\multiply\height by 2
\putvmorphism(\xpos,\ypos)[#4`#6`{#7}]{\height}{\arrowtypea}l}%
}}

\def\putDtriangle{\@ifnextchar <{\putDtrianglep}{\putDtrianglep
    <\arrowtypea`\arrowtypeb`\arrowtypec;\height>}}
\def\Dtriangle{\@ifnextchar <{\Dtrianglep}{\Dtrianglep
   <\arrowtypea`\arrowtypeb`\arrowtypec;\height>}}
\def\Dtrianglep<#1>[#2`#3`#4;#5`#6`#7]{{%%  %          | \
\settriparms[#1]%                           %          |  \#6
\width=\height                              %          |   \
\diagram%                                   %          |    v
\putDtrianglep<\arrowtypea`\arrowtypeb`%    %        #5|    #3
\arrowtypec;\height>%                       %          |    /
(0,0)[#2`#3`#4;#5`#6`{#7}]%                 %          |   /#7
\enddiagram%%                               %          |  /
}}                                          %          v v
\def\setrecparms[#1`#2]{\width=#1 \height=#2}%
\def\recursep<#1`#2>[#3;#4`#5`#6`#7`#8]{{%
\width=#1 \height=#2
\settokens[#3]
\settowidth{\tempdimen}{$\tokena$}
\ifdim\tempdimen=0pt
  \savebox{\tempboxa}{\hbox{$\tokenb$}}%
  \savebox{\tempboxb}{\hbox{$\tokend$}}%
  \savebox{\tempboxc}{\hbox{$#6$}}%
\else
  \savebox{\tempboxa}{\hbox{$\hbox{$\tokena$}\times\hbox{$\tokenb$}$}}%
  \savebox{\tempboxb}{\hbox{$\hbox{$\tokena$}\times\hbox{$\tokend$}$}}%
  \savebox{\tempboxc}{\hbox{$\hbox{$\tokena$}\times\hbox{$#6$}$}}%
\fi
\ypos=\height
\divide\ypos by 2
\xpos=\ypos
\advance\xpos by \width
\bfig
\putCtrianglep<-1`1`1;\ypos>(0,0)[`\tokenc`;#5`#6`{#7}]%
\puthmorphism(\ypos,0)[\tokend`\usebox{\tempboxb}`{#8}]{\width}{-1}b%
\puthmorphism(\ypos,\height)[\tokenb`\usebox{\tempboxa}`{#4}]{\width}{-1}a%
\advance\ypos by \width
\putvmorphism(\ypos,\height)[``\usebox{\tempboxc}]{\height}1r%
\efig
}}

\def\recurse{\@ifnextchar <{\recursep}{\recursep<\width`\height>}}

\def\puttwohmorphisms(#1,#2)[#3`#4;#5`#6]#7#8#9{{%
\puthmorphism(#1,#2)[#3`#4`]{#7}0a
\ypos=#2
\advance\ypos by 20
\puthmorphism(#1,\ypos)[\phantom{#3}`\phantom{#4}`#5]{#7}{#8}a
\advance\ypos by -40
\puthmorphism(#1,\ypos)[\phantom{#3}`\phantom{#4}`#6]{#7}{#9}b
}}

\def\puttwovmorphisms(#1,#2)[#3`#4;#5`#6]#7#8#9{{%
\putvmorphism(#1,#2)[#3`#4`]{#7}0a
\xpos=#1
\advance\xpos by -20
\putvmorphism(\xpos,#2)[\phantom{#3}`\phantom{#4}`#5]{#7}{#8}l
\advance\xpos by 40
\putvmorphism(\xpos,#2)[\phantom{#3}`\phantom{#4}`#6]{#7}{#9}r
}}

\def\puthcoequalizer(#1)[#2`#3`#4;#5`#6`#7]#8#9{{%
\setpos(#1)%
\puttwohmorphisms(\xpos,\ypos)[#2`#3;#5`#6]{#8}11%
\advance\xpos by #8
\puthmorphism(\xpos,\ypos)[\phantom{#3}`#4`#7]{#8}1{#9}
}}

\def\putvcoequalizer(#1)[#2`#3`#4;#5`#6`#7]#8#9{{%
\setpos(#1)%
\puttwovmorphisms(\xpos,\ypos)[#2`#3;#5`#6]{#8}11%
\advance\ypos by -#8
\putvmorphism(\xpos,\ypos)[\phantom{#3}`#4`#7]{#8}1{#9}
}}

\def\putthreehmorphisms(#1)[#2`#3;#4`#5`#6]#7(#8)#9{{%
\setpos(#1) \settypes(#8)
\if a#9 %
     \vertsize{\tempcounta}{#5}%
     \vertsize{\tempcountb}{#6}%
     \ifnum \tempcounta<\tempcountb \tempcounta=\tempcountb \fi
\else
     \vertsize{\tempcounta}{#4}%
     \vertsize{\tempcountb}{#5}%
     \ifnum \tempcounta<\tempcountb \tempcounta=\tempcountb \fi
\fi
\advance \tempcounta by 60
\puthmorphism(\xpos,\ypos)[#2`#3`#5]{#7}{\arrowtypeb}{#9}
\advance\ypos by \tempcounta
\puthmorphism(\xpos,\ypos)[\phantom{#2}`\phantom{#3}`#4]{#7}{\arrowtypea}{#9}
\advance\ypos by -\tempcounta \advance\ypos by -\tempcounta
\puthmorphism(\xpos,\ypos)[\phantom{#2}`\phantom{#3}`#6]{#7}{\arrowtypec}{#9}
}}

\def\setarrowtoks[#1`#2`#3`#4`#5`#6]{%
\def\toka{#1}
\def\tokb{#2}
\def\tokc{#3}
\def\tokd{#4}
\def\toke{#5}
\def\tokf{#6}
}
\def\hex{\@ifnextchar <{\hexp}{\hexp<1000`400>}}
\def\hexp<#1`#2>[#3`#4`#5`#6`#7`#8;#9]{%
\setarrowtoks[#9]
\yext=#2 \advance \yext by #2
\xext=#1 \advance\xext by \yext
\bfig
\putCtriangle<-1`0`1;#2>(0,0)[`#5`;\tokb``\tokd]
\xext=#1 \yext=#2 \advance \yext by #2
\putsquare<1`0`0`1;\xext`\yext>(#2,0)[#3`#4`#7`#8;\toka```\tokf]
\advance \xext by #2
\putDtriangle<0`1`-1;#2>(\xext,0)[`#6`;`\tokc`\toke]
\efig
}

\QQQ{Language}{
American English
}

\begin{document}

\title{{\bf \ Brown-Peterson spectra}\\
{\bf \ in stable }${\Bbb A}^1${\bf -homotopy theory}}
\author{Gabriele Vezzosi\thanks{%
This paper was written during a very pleasant visit at the Institute for
Advanced Study, Princeton, NJ. The author wishes to thank all people there
for creating such an exceptional atmosphere.} \\
%EndAName
Dipartimento di Matematica, Universit\`a di Bologna,\\
Piazza di Porta S. Donato 5, Bologna 40127 - Italy}
\maketitle

\begin{abstract}
We characterize ring spectra morphisms from the algebraic cobordism spectrum 
${\Bbb MGL}$ (\cite{Vo1}) to an oriented spectrum ${\Bbb E}$ (in the sense
of Morel \cite{Mo}) via formal group laws on the ''topological'' subring $%
E^{*}=\oplus _iE^{2i,i}$ of $E^{**}$. This result is then used to construct
for any prime $p$ a motivic Quillen idempotent on ${\Bbb MGL}_{(p)}$. This
defines the $BP$-spectrum associated to the prime $p$ as in Quillen's \cite
{Q1} for the complex-oriented topological case.

\smallskip\ 

\end{abstract}

\section{Introduction}

My interest in the subject of this paper originated from the idea to extend
Totaro's construction (\cite{To}) of a refined cycle map with values in a
quotient of complex cobordism, into Voevodsky's algebraic cobordism setting (%
\cite{Vo1}). To construct such an extension, one need to prove, as Totaro
did for complex cobordism, a slight refinement of Quillen theorem 5.1 of 
\cite{Q2}. As Quillen's original proof rests on some finiteness results
(Prop. 1.12, \cite{Q2}) together with a geometric interpretation of complex
cobordism groups and since these conditions are not currently\footnote{%
Actually, the first is genuinely false.} available in algebraic cobordism,
we decided to follow Totaro's argument which uses $BP$ and $BP\left\langle
n\right\rangle $ spectra.

The first step, i.e. to define $BP$-spectra in ${\Bbb A}^1$-homotopy theory,
is carried over in this paper. The other steps (construction of $%
BP\left\langle n\right\rangle $ spectra and Wilson's Theorem, \cite{W}, p.
413) will be considered in subsequent papers.

The intuition guiding the construction of $BP$-spectra obviously comes from
topology although some extra care is required in ${\Bbb A}^1$-homotopy
theory because of the lacking of some facts which are well known in
topology. I am referring in particular to the fact that ${\rm MGL}^{**}={\rm %
MGL}^{**}\left( {\rm Spec}\left( k\right) \right) $ has not been computed
yet, although what I call in this paper its ''topological''subring ${\rm MGL}%
^{*}\doteq \oplus _i{\rm MGL}^{2i,i}$, is conjectured to be isomorphic to
the complex cobordism ring ${\rm MU}^{*}$ and therefore, by Quillen theorem,
to the Lazard ring (\cite{Vo1} 6.3).

Here is a brief description of the contents. The first section is devoted to
explain the algebraic geometric motivation that lead to the problem of
defining Brown-Peterson spectra in ${\Bbb A}^1$-homotopy theory. In the
second section we use some results of Morel (\cite{Mo}), to establish the
basic properties of oriented ring spectra in the stable ${\Bbb A}^1$%
-homotopy category ${\cal SH}\left( k\right) $ of smooth schemes over a
ground field $k$ (\cite{Vo1}, \cite{MV}), especially the Thom isomorphism in
its stable and unstable versions. In the third section, for an oriented
spectrum ${\Bbb E}$, we prove the equivalence between giving an orientation
on ${\Bbb E}$, giving a ring spectra map ${\Bbb MGL\rightarrow E}$ and
giving a formal group law on ${\rm E}^{*}=\oplus _i{\rm E}^{2i,i}\left( {\rm %
Spec}\left( k\right) \right) $ isomorphic to the one associated to the given
orientation. This enables us to give the construction of $BP$-spectra in the
fourth section through the construction of a ''motivic'' Quillen idempotent.

We remark that all our proofs work also in the topological case, with
algebraic cobordism replaced by complex cobordism and that the safe path we
have to follow in order to avoid any use of Quillen theorem on the Lazard
ring, yields proofs which in the topological case are different from the
standard ones.

This paper was ''triggered'' by a question Vladimir Voevodsky posed to me
during a visit at the Instistute for Advanced Study, in March 2000. I wish
to thank him for helpful discussions and advises on this and other topics. I
am also indebted to Fabien Morel for allowing me to use the results in \cite
{Mo}.

It is also a pleasure to thank Dan Christensen, Charles Rezk, Burt Totaro
and Chuck Weibel for listening carefully and patiently to a non-topologist
speaking about topology and for their useful advises.

\section{A motivation}

In this section we briefly explain the motivations that lead us to the
problem of finding a definition of Brown-Peterson spectra in ${\Bbb A}^1$%
-homotopy theory. This section is mainly conjectural, independent from the
others and only meant to suggest one possible way to approach the problems
described below. 

In \cite{To}, Totaro defined a refined cycle class map 
\begin{equation}
A_{*}\left( X\right) \stackrel{\widetilde{{\rm cl}}_X}{\longrightarrow }%
\widetilde{{\rm MU}}_{2*}\left( X^{{\rm an}}\right) \doteq \left( {\rm MU}%
_{*}\left( X^{{\rm an}}\right) \otimes _{^{{\rm MU}_{*}}}{\Bbb Z}\right)
_{2*}  \label{tot}
\end{equation}
where $X$ is an equidimensional smooth complex algebraic scheme, $X^{{\rm an}%
}$ the associated complex manifold, $A_{*}\left( X\right) $ is the Chow
group of $X$ (\cite{Fu}, Ch. 1), ${\rm MU}_{*}\left( Y\right) $ denotes the
complex bordism group of the topological space $Y$ and 
\[
{\rm MU}_{*}={\rm MU}_{*}\left( {\rm pt}\right) \simeq {\Bbb Z}\left[
x_1,x_{2,\ldots ,}x_{n,\ldots }\right] 
\]
with $\deg x_i=2i$. Here ${\Bbb Z}$ is considered as an ${\rm MU}_{*}$%
-module via the map sending each $x_i$ to $0$ i.e. 
\[
\widetilde{{\rm MU}}_{*}\left( X^{{\rm an}}\right) ={\rm MU}_{*}\left( X^{%
{\rm an}}\right) \otimes _{^{{\rm MU}_{*}}}{\Bbb Z\simeq }\frac{{\rm MU}%
_{*}\left( X^{{\rm an}}\right) }{{\rm MU}^{>0}\cdot {\rm MU}_{*}\left( X^{%
{\rm an}}\right) ^{}}\text{.}
\]
Totaro proved that the classical cycle map ${\rm cl}_X:A_{*}\left( X\right)
\rightarrow {\rm H}_{2*}^{{\rm BM}}\left( X^{{\rm an}};{\Bbb Z}\right) $ (%
\cite{Fu}, 19.1) to Borel-Moore homology factors as 
\[
A_{*}\left( X\right) \stackrel{\widetilde{{\rm cl}}_X}{\longrightarrow }%
\widetilde{{\rm MU}}_{2*}\left( X^{{\rm an}}\right) \stackrel{\gamma _{X^{%
{\rm an}}}}{\longrightarrow }{\rm H}_{2*}^{{\rm BM}}\left( X^{{\rm an}};%
{\Bbb Z}\right) 
\]
where $\gamma _{X^{{\rm an}}}$ is induced by the canonical map of homology
theories 
\[
{\rm MU}_{*}\left( X^{{\rm an}}\right) \longrightarrow {\rm H}_{*}^{{\rm BM}%
}\left( X^{{\rm an}};{\Bbb Z}\right) 
\]
\[
\left[ f:M\rightarrow X^{{\rm an}}\right] \longmapsto f_{*}\left( \eta
_M\right) 
\]
where we used the geometric interpretation of complex bordism classes as
equivalence classes of proper maps from weakly complex manifolds (\cite{St})
and $\eta _M$ denotes the fundamental class of the weakly complex manifold $M
$.

The refined cycle map $\widetilde{{\rm cl}}_X$ sends the class of a cycle $%
Z\hookrightarrow X$ to the class of the composition 
\[
\widetilde{Z}\longrightarrow Z\hookrightarrow X
\]
$\widetilde{Z}\rightarrow Z$ being any resolution of singularities. The
reason it is well defined is a combination of  Hironaka's theorem,
Poincar\'e duality and the following

\begin{theorem}
\label{TQ}{\rm (Quillen-Totaro theorem, \cite{To} Thm. 2.2)}

Let $Y$ be a finite cell complex. Then the canonical map 
\[
\widetilde{{\rm MU}}\Sp  \\ \ast  \endSp \left( Y\right) \stackrel{\gamma _Y%
}{\longrightarrow }H^{*}\left( Y;{\Bbb Z}\right) 
\]
is injective in degrees $\leq 2$.
\end{theorem}

In degrees $\leq 0$ this is a consequence of Quillen theorem (\cite{Q2},
5.1) but Totaro's proofs uses a different approach, through $BP$-spectra $%
{\Bbb BP}$, truncated $BP$-spectra ${\Bbb BP}\left\langle n\right\rangle $
and Wilson theorem (\cite{W}, p. 118).

Our question is whether it exists a generalization of Totaro's refined cycle
map over an arbitrary field $k$ admitting resolution of singularities., with
complex (co)bordism replaced by algebraic (co)bordism. It turns out that
even the formulation of the analog of Quillen theorem requires a little care.

Following \cite{Vo1}, let us denote by ${\Bbb H}_{{\Bbb Z}}$ the motivic
Eilenberg-Mac Lane spectrum and by ${\Bbb MGL\doteq }\left( {\rm MGL}\left(
n\right) \right) _n$ the algebraic cobordism spectrum (\cite{Vo1}, 6.1 and
6.3). These are objects in the stable ${\Bbb A}^1$-homotopy category ${\cal %
SH}\left( k\right) $ of smooth schemes over $k$. First of all, by \cite{Vo3}
Thm. 3.21
\[
{\rm Hom}_{{\cal SH}\left( k\right) }\left( {\Bbb MGL}{\bf ,}{\Bbb H}_{{\Bbb %
Z}}\right) \doteq {\rm H}_{{\Bbb Z}}^{0,0}\left( {\Bbb MGL}\right) \simeq 
{\Bbb Z}
\]
canonically and if we denote by $\tau $ a generator of ${\rm H}_{{\Bbb Z}%
}^{0,0}\left( {\Bbb MGL}\right) $, for any smooth scheme $X$ over $k$ we
have an induced morphism of cohomology theories
\[
\tau _X:{\rm MGL}^{**}\left( X\right) \longrightarrow {\rm H}_{{\Bbb Z}%
}^{**}\left( X\right) .
\]
Let 
\[
{\rm MGL}^{*}\left( X\right) \doteq \oplus _i{\rm MGL}^{2i,i}\left( X\right) 
\]
and
\[
{\rm H}_{{\Bbb Z}}^{*}\left( X\right) \doteq \oplus _i{\rm H}_{{\Bbb Z}%
}^{2i,i}\left( X\right) .
\]
Since (\cite{We}, p. 293) ${\rm H}_{{\Bbb Z}}^p\doteq {\rm H}_{{\Bbb Z}%
}^p\left( {\rm Spec}k\right) \simeq A^p\left( {\rm Spec}k\right) =0$ if $p<0$%
, the restriction of $\tau _X$ to ${\rm MGL}^{*}\left( X\right) $ factors
through
\begin{equation}
\widetilde{\tau _X}:\widetilde{{\rm MGL}}^{*}\left( X\right) \doteq \frac{%
{\rm MGL}^{*}\left( X\right) }{{\rm MGL}^{<0}\cdot {\rm MGL}^{*}\left(
X\right) }\longrightarrow {\rm H}_{{\Bbb Z}}^{*}\left( X\right) 
\label{cippa}
\end{equation}
as in the complex (oriented) case.

\begin{remark}
{\rm \ The ring }${\rm MGL}^{**}={\rm MGL}^{**}\left( {\rm Spec}k\right) $%
{\rm \ is not known but it is conjectured that its subring }${\rm MGL}^{*}$%
{\rm \ is isomorphic to }${\rm MU}^{*}${\rm \ (\cite{Vo1}, 6.3). A part of
this conjecture, i.e. that }${\rm MGL}^{*}${\rm \ is zero in positive
degrees, follows immediately from the }${\cal SH}\left( k\right) ${\rm %
-version of the  Connectivity Theorem 4.14 in \cite{Vo1}. The rest of the
conjecture is, as far as we know, still open. This is the main reason why we
will need to avoid the relation of  }${\rm MGL}^{*}${\rm \ to the Lazard
ring in the following sections.} 
\end{remark}

Now, if $\eta _Y\in {\rm MGL}_{2n,n}\left( Y\right) $ denotes the
fundamental class of an $n$-dimensional smooth $k$-scheme $Y$ in algebraic
cobordism (\cite{Vo3}), the motivic analog of Totaro's refined cycle map
should be the map
\begin{equation}
\widetilde{{\rm CL}}_X:A_{*}\left( X\right) \longrightarrow \widetilde{{\rm %
MGL}}_{*}\left( X\right) =\frac{{\rm MGL}_{*}\left( X\right) }{{\rm MGL}%
_{>0}\cdot {\rm MGL}_{*}\left( X\right) }  \label{mot}
\end{equation}
sending the class $\left[ Z\hookrightarrow X\right] $ of a cycle of
dimension $i$ in $X$ to the class modulo ${\rm MGL}_{>0}\cdot {\rm MGL}%
_{*}\left( X\right) $ of 
\[
\left( \left( {\Bbb P}^1,\infty \right) ^{\wedge i}\stackrel{\eta _{%
\widetilde{Z}}}{\longrightarrow }\Sigma ^\infty \widetilde{Z}_{+}\wedge 
{\Bbb MGL}\stackrel{\Sigma ^\infty f\wedge {\rm id}}{\longrightarrow }\Sigma
^\infty X_{+}\wedge {\Bbb MGL}\right) \in {\rm MGL}_{2i,i}\left( X\right) 
\]
where 
\[
f:\widetilde{Z}\rightarrow Z\hookrightarrow X
\]
$\widetilde{Z}\rightarrow Z$ being any resolution of singularities.
Admitting a Poincar\'e duality for algebraic (co)bordism and motivic
(co)homology, well definiteness of (\ref{mot}) would be a consequence of the
following 

\begin{conjecture}
{\rm (Motivic Quillen-Totaro Theorem)}

For any smooth scheme over $k$, the map {\rm (\ref{cippa})}
\[
\widetilde{\tau _X}:\widetilde{{\rm MGL}}^{*}\left( X\right) \longrightarrow 
{\rm H}_{{\Bbb Z}}^{*}\left( X\right) 
\]
is injective in degrees $\leq 2$.
\end{conjecture}

Note that in this case, the ''classical'' cycle map should be the identity
\[
{\rm CL}_X={\rm id}:A_{*}\left( X\right) \rightarrow {\rm H}_{*}^{{\Bbb Z}%
}\left( X\right) \simeq A_{*}\left( X\right) 
\]
(the last isomorphism given by Poincar\'e duality) and the factorization $%
{\rm id}=\widetilde{\tau ^X}\circ \widetilde{{\rm CL}}_X$ would imply that $%
\widetilde{\tau ^X}$ is surjective and $\widetilde{{\rm CL}}_X$ injective.

Following the idea in Totaro's proof of \cite{To}, Thm. 2.2, the first thing
to know is how to construct $BP$-spectra in ${\cal SH}\left( k\right) $.
This is done in the following sections.

\section{Oriented spectra and Thom isomorphism}

Throughout the paper we fix a base field $k$ and work in the stable ${\Bbb A}%
^1$-homotopy category ${\cal SH}\left( k\right) $ of smooth schemes over $k$%
, as described in \cite{Vo1}, whose notations we follow closely. In
particular, $\Sigma ^\infty $ will denote the infinite $\left( {\Bbb P}%
^1,\infty \right) $-suspension, and for any space $X$ over $k$, we write $%
X_{+}$ for the space $X\coprod {\rm Spec}\left( k\right) $ pointed by ${\rm %
Spec}\left( k\right) $. $S^0$ will denote ${\rm Spec}\left( k\right) _{+}$
and ${\Bbb S}^{p,q}$ will denote the spectrum obtained from the smash
product of the mixed spheres $S_s^{p-q}$ and $S_t^q$ as described in \cite
{MV} 3.2.2; recall that for any spectrum ${\Bbb E}$ in ${\cal SH}\left(
k\right) $ and any space $X$ (respectively, spectrum ${\Bbb F}$) we have a
cohomology theory 
\[
E^{p,q}\left( X\right) \doteq {\rm Hom}_{{\cal SH}\left( k\right) }\left(
\Sigma ^\infty \left( X_{+}\right) ,{\Bbb S}^{p,q}\wedge {\Bbb E}\right) 
\]
(respectively, 
\[
E^{p,q}\left( {\Bbb F}\right) \doteq {\rm Hom}_{{\cal SH}\left( k\right)
}\left( {\Bbb F},{\Bbb S}^{p,q}\wedge {\Bbb E}\right) \text{ ).} 
\]
When no reasonable ambiguity seems to take place, we also write simply $%
{\Bbb P}^1$ for the pointed space $\left( {\Bbb P}^1,\infty \right) $.

In this section we mainly follows \cite{Mo} and draw some consequences
thereof.

If ${\Bbb E}$ is a ring spectrum in ${\cal SH}\left( k\right) $ (in the weak
sense), there is a canonical element $x_E^0\in E^{2,1}\left( {\Bbb P}%
^1\right) $ given by the composition 
\[
\Sigma ^\infty {\Bbb P}_{+}^1\longrightarrow \Sigma ^\infty \left( {\Bbb P}%
^1,\infty \right) \widetilde{\longrightarrow }\Sigma ^\infty \left( {\Bbb P}%
^1,\infty \right) \wedge S^0\stackrel{id\wedge \eta }{\longrightarrow }%
\Sigma ^\infty \left( {\Bbb P}^1,\infty \right) \wedge {\Bbb E} 
\]
where $\eta :\Sigma ^\infty \left( {\Bbb P}^1,\infty \right) \rightarrow 
{\Bbb E}$ denotes the unit morphism of the ring spectrum ${\Bbb E}$.

\begin{definition}
{\rm (\cite{Mo}, 3.2.3)}

An {\rm oriented} ring spectrum in ${\cal SH}\left( k\right) $ is a pair $%
\left( {\Bbb E},x_E\right) $ where ${\Bbb E}$ is a commutative ring spectrum
and $x_E$ is an element in $E^{2,1}\left( {\Bbb P}^\infty \right) $
restricting to the canonical element $x_E^0$ along the canonical inclusion $%
{\Bbb P}^1\rightarrow {\Bbb P}^\infty $.
\end{definition}

Here, ${\Bbb P}^\infty \doteq $ {\rm colim}$_n{\Bbb P}^n$.

If ${\rm Gr}_{n,N}$ denote the Grassmannian of $n$-planes in ${\Bbb A}^N$, $%
N>n$, let us denote by ${\rm BGL}_n$ (respectively, {\rm BGL}) the infinite
Grassmannian {\rm colim}$_N\left( {\rm Gr}_{n,N}\right) $ of $n$-planes
(respectively, {\rm colim}$_n\left( {\rm BGL}_n\right) $). Moreover, let us
denote by ${\Bbb MGL\doteq }\left( {\rm MGL}\left( n\right) \right) _n$ the
algebraic cobordism spectrum (\cite{Vo1}, 6.3).

\begin{lemma}
\label{MGL1}The zero section map 
\[
s_0:{\Bbb P}^\infty ={\rm BGL}_1\longrightarrow {\rm MGL}\left( 1\right) 
\]
is a weak equivalence.
\end{lemma}

\TeXButton{Proof}{\proof} For any $n>0$, the closed immersion ${\Bbb P}%
^{n-1}\hookrightarrow {\Bbb P}^n$ has normal bundle the canonical line $%
{\cal L}_{n-1}$ bundle on ${\Bbb P}^{n-1}$ and ${\Bbb P}^n-{\Bbb P}^{n-1}$
is isomorphic to ${\Bbb A}^n$; hence (\cite{MV}, Th. 3.2.23) the Thom space $%
{\rm Th}\left( {\cal L}_{n-1}\right) $ (\cite{Vo1}, p. 422) is weakly
equivalent to ${\Bbb P}^n$ and these weak equivalences are compatible with
respect to the maps in the direct system $\left\{ \cdots \hookrightarrow 
{\Bbb P}^{n-1}\hookrightarrow {\Bbb P}^n\hookrightarrow {\Bbb P}%
^{n+1}\hookrightarrow \cdots \right\} $. The result follows by passing to
the colimit. \TeXButton{End Proof}{\endproof}

\begin{remark}
\label{canonical}{\rm Note that the algebraic cobordism spectrum }${\Bbb MGL}
${\rm \ has a canonical orientation }$x_{{\Bbb MGL}}${\rm \ given by the
composition } 
\[
\Sigma ^\infty {\Bbb P}_{+}^\infty \longrightarrow \Sigma ^\infty {\Bbb P}%
^\infty \stackrel{\Sigma ^\infty \left( s_0\right) }{\longrightarrow }\Sigma
^\infty {\rm MGL}\left( 1\right) \stackrel{\nu }{\longrightarrow }\Sigma
^\infty {\Bbb P}^1\wedge {\Bbb MGL}
\]
{\rm where }$\nu ${\rm \ is defined using the bonding maps of the }${\Bbb MGL%
}${\rm -spectrum as } 
\[
\left( {\Bbb P}^1\right) ^{\wedge n}\wedge {\rm MGL}\left( 1\right)
\rightarrow \left( {\Bbb P}^1\right) ^{\wedge n-1}\wedge {\rm MGL}\left(
2\right) \rightarrow \cdots \rightarrow {\Bbb P}^1\wedge {\rm MGL}\left(
n\right) \text{.}
\]
\end{remark}

\begin{proposition}
\label{morel}{\rm (\cite{Mo}, 3.2.9)}

{\rm (i)} If $\left( {\Bbb E},x_E\right) $ is an oriented ring spectrum, for
any space $X$ over $k$, we have 
\[
\alpha \cup \beta =\left( -1\right) ^{pp^{\prime }}\beta \cup \alpha 
\]
for any $\alpha \in E^{p,q}\left( X\right) $ and $\beta \in E^{p^{\prime
},q^{\prime }}\left( X\right) $. In particular, the subring $E^{*}\doteq
\oplus _iE^{2i,i}\left( {\rm Spec}\left( k\right) \right) $ is commutative.

{\rm (ii) }For any $n>0$, there is a canonical isomorphism 
\[
E^{**}\left( {\Bbb P}^n\right) \simeq E^{**}\left[ x_E\right] /\left(
x_E^{n+1}\right) 
\]
where we still denote by $x_E$ its pullback along ${\Bbb P}^n\hookrightarrow 
{\Bbb P}^\infty $.

{\rm (iii)} Moreover, for any $n>0$, there is a canonical isomorphism 
\[
E^{**}\left( \stackunder{n\text{ times}}{\underbrace{{\Bbb P}^\infty \times
\cdots \times {\Bbb P}^\infty }}\right) \simeq E^{**}\left[ \left[
x_1,\ldots x_n\right] \right] 
\]
where $x_i$ denote the pullback of $x_E$ along the $i$-th projection ${\rm pr%
}_i:{\Bbb P}^\infty \times \cdots \times {\Bbb P}^\infty \rightarrow {\Bbb P}%
^\infty $, $i=1,\ldots ,n$.
\end{proposition}

Note that (iii), follows from (ii) (which is \cite{Mo}, 3.2.9, (2)), by
passing to the limit after recognizing the Mittag-Leffler condition is
satisfied.

\begin{proposition}
\label{lullo}If $\left( {\Bbb E},x_E\right) $ is an oriented ring spectrum,
for any $n>0$, the pullback along the canonical map 
\begin{equation}
\stackunder{n\text{ times}}{\underbrace{{\Bbb P}^\infty \times \cdots \times 
{\Bbb P}^\infty }}\longrightarrow {\rm BGL}_n  \label{bgln}
\end{equation}
classifying the product of canonical line bundles {\rm (}respectively, the
pullback along the colimit{\rm \footnote{{\rm The colimit over the inclusion
of any rational point in ${\Bbb P}^\infty $. Any choice will yield the same
direct system up to weak equivalences since any two such points belong to an
affine line over $k$.}}} of maps {\rm (\ref{bgln})} 
\[
\theta :\left( {\Bbb P}^\infty \right) ^\infty \doteq \text{{\rm colim}}%
_n\left( {\Bbb P}^\infty \right) \longrightarrow {\rm BGL}\text{{\rm )}}
\]
induces a monomorphism 
\[
E^{**}\left( {\rm BGL}_n\right) \simeq E^{**}\left[ \left[ c_1,\ldots
,c_n\right] \right] \hookrightarrow E^{**}\left[ \left[ x_1,\ldots
,x_n\right] \right] \simeq E^{**}\left( \stackunder{n\text{ times}}{%
\underbrace{{\Bbb P}^\infty \times \cdots \times {\Bbb P}^\infty }}\right) 
\]
{\rm (}respectively, 
\[
E^{**}\left( {\rm BGL}\right) \simeq E^{**}\left[ \left[ c_1,\ldots
,c_n,\ldots \right] \right] \hookrightarrow E^{**}\left[ \left[ x_1,\ldots
,x_n,\ldots \right] \right] \simeq E^{**}\left( \left( {\Bbb P}^\infty
\right) ^\infty \right) \text{{\rm )}}
\]
where $c_i$ denotes the $i$-th elementary symmetric function on the $x_j$'s.
\end{proposition}

\TeXButton{Proof}{\proof} Both the asserts follows from \cite{Mo} 3.2.10
(2), by induction on $n$. \TeXButton{End Proof}{\endproof}

\begin{corollary}
\label{thomclasses}If $\left( {\Bbb E},x_E\right) $ is an oriented ring
spectrum, there exists a canonical ''family of universal Thom classes'' {\rm %
(\cite{CGK})} i.e. a family $\left( \tau _n^E\right) _{n>0}$ with $\tau
_n^E\in E^{2n,n}\left( {\rm MGL}\left( n\right) \right) $, such that:

{\rm (i)} $\tau _1^E=x_E$ {\rm (}this makes sense because of Lemma {\rm \ref
{MGL1})};

{\rm (ii)} the family is multiplicative in the sense that $\tau _{n+m}^E$
pulls back to $\tau _n^E\wedge \tau _m^E$ along the map ${\rm MGL}\left(
n\right) \wedge {\rm MGL}\left( m\right) \rightarrow {\rm MGL}\left(
n+m\right) $ (induced by the canonical map ${\rm BGL}_n\times {\rm BGL}%
_m\rightarrow {\rm BGL}_{n+m}$making ${\rm BGL}$ into an $H$-space).
\end{corollary}

\TeXButton{Proof}{\proof}If $\xi _n\rightarrow {\rm BGL}_n$ denotes the
universal $n$-plane bundle, we have a cofiber sequence (\cite{MV}, 3.2.17) 
\[
{\Bbb P}\left( \xi _n\right) \longrightarrow {\Bbb P}\left( \xi _n\oplus 
{\bf 1}\right) \longrightarrow {\rm MGL}\left( n\right) 
\]
whose associated $E$-cohomology long exact sequence yields, by the
projective bundle theorem (\cite{Mo}, 3.2.10, (1)), a short exact sequence 
\[
0\rightarrow E^{**}\left( {\rm MGL}\left( n\right) \right) \longrightarrow
E^{**}\left( {\rm BGL}_n\right) \left[ t\right] /\left( t\cdot f\left(
t\right) \right) \stackrel{\pi }{\longrightarrow }E^{**}\left( {\rm BGL}%
_n\right) \left[ s\right] /\left( f\left( s\right) \right) \rightarrow 0 
\]
where $f\left( t\right) \doteq t^n+c_1\left( \xi _n\right) t^{n-1}+\cdots
+c_n\left( \xi _n\right) $, $t$ corresponding to the first Chern class of
the canonical line bundle on ${\Bbb P}\left( \xi _n\right) $, ${\bf 1}$
denotes the trivial line bundle over ${\rm BGL}_n$, $\pi $ maps $t$ to $s$
and $c_i\left( {\cal E}\right) \in E^{2i,i}\left( {\rm BGL}_n\right) $
denotes here the $i$-th Chern class of a vector bundle ${\cal E}$ (we have
used that $c_i\left( \xi _n\oplus 1\right) =c_i\left( \xi _n\right) $ . Then 
$f\left( t\right) $ has bidegree $\left( 2n,n\right) $ and is in the kernel
of $\pi $. Then, $\tau _n^E$ is the unique element in $E^{2n,n}\left( {\rm %
MGL}\left( n\right) \right) $ mapping to $f\left( t\right) $ and it is easy
to verify that the family $\left( \tau _n^E\right) _{n>0}$ defined in this
way is multiplicative. \TeXButton{End Proof}{\endproof}

\begin{remark}
{\rm The universal Thom classes }$\left( \tau _n^E\right) _{n>0}${\rm \
admit the following equivalent characterization For any }$n>0${\rm ,
consider the product }$c_n=x_1\cdots x_n\in E^{2n,n}\left( {\rm BGL}%
_n\right) ${\rm \ (Prop. \ref{BG}). By Prop. \ref{BG}, the }$E${\rm %
-cohomology long exact sequence associated to the cofiber sequence } 
\[
{\rm BGL}_{n-1}\longrightarrow {\rm BGL}_n\longrightarrow {\rm MGL}\left(
n\right) 
\]
{\rm yields a short exact sequence } 
\[
0\rightarrow E^{2n,n}\left( {\rm MGL}\left( n\right) \right) \stackrel{%
\varphi }{\longrightarrow }E^{2n,n}\left( {\rm BGL}_n\right) \stackrel{\psi 
}{\longrightarrow }E^{2n,n}\left( {\rm BGL}_{n-1}\right) \rightarrow 0
\]
{\rm where }$\psi ${\rm \ maps }$g\left( c_{1,\ldots ,}c_n\right) ${\rm \ to 
}$g\left( c_{1,\ldots ,}c_{n-1},0\right) ${\rm . Therefore there exists a
unique class } 
\[
\tau _n^E\in E^{2n,n}\left( {\rm MGL}\left( n\right) \right) 
\]
{\rm such that }$\varphi \left( \tau _n^E\right) =c_n${\rm .}
\end{remark}

\smallskip\ 

As in the topological case, the projective bundle structure theorem for $E$%
-cohomology (\cite{Mo}, 3.2.10, (1)) implies the Thom isomorphism.

Let $\left( {\Bbb E},x_E\right) $ be an oriented ring spectrum, $X$ a smooth
scheme over $k$ and ${\cal E}\rightarrow X$ be a vector bundle of rank $r$.
If ${\rm Th}\left( {\cal E}/X\right) $ denotes the Thom space of ${\cal E}$ (%
\cite{Vo1}, p. 422), the diagonal map $\delta :X\rightarrow X\times X$
induces a Thom diagonal 
\begin{equation}
\Delta _{{\cal E}}:{\rm Th}\left( {\cal E}/X\right) \longrightarrow {\rm Th}%
\left( {\cal E}/X\right) \wedge X_{+}\text{ .}  \label{diagonal}
\end{equation}
Since ${\cal E}$ has rank $r$, there is a canonical map 
\[
\lambda _{{\cal E}}:{\rm Th}\left( {\cal E}/X\right) \longrightarrow {\rm MGL%
}\left( r\right) . 
\]
Therefore we have a Thom map 
\begin{equation}
\Phi _{{\cal E}}:E^{**}\left( X\right) \longrightarrow E^{*\text{ }+\text{ }%
2r,\,*\text{ }+\text{ }r}\left( {\rm Th}\left( {\cal E}/X\right) \right)
\label{thomiso}
\end{equation}
which assigns to an element $\alpha \in E^{p,q}\left( X\right) ={\rm Hom}_{%
{\cal SH}\left( k\right) }\left( \Sigma ^\infty \left( X_{+}\right) ,{\Bbb S}%
^{p,q}\wedge {\Bbb E}\right) $ the element $\Phi _{{\cal E}}\left( \alpha
\right) $ in $E^{p+2r,\,q+r}\left( {\rm Th}\left( {\cal E}/X\right) \right) $
given by the composition 
\[
\Sigma ^\infty {\rm Th}\left( {\cal E}/X\right) \stackrel{\Sigma ^\infty
\Delta _{{\cal E}}}{\longrightarrow }\Sigma ^\infty \left( {\rm Th}\left( 
{\cal E}/X\right) \wedge X_{+}\right) \stackrel{\Sigma ^\infty \left(
\lambda _{{\cal E}}\wedge \alpha \right) }{\longrightarrow }\Sigma ^\infty
\left( {\rm MGL}\left( r\right) \wedge {\Bbb S}^{p,q}\wedge {\Bbb E}\right)
\rightarrow 
\]
\[
\stackrel{\tau _r^E\wedge id}{\longrightarrow }{\Bbb S}^{2r,r}\wedge {\Bbb %
E\wedge S}^{p,q}\wedge {\Bbb E\longrightarrow S}^{p+2r,q+r}\wedge {\Bbb E} 
\]
where $\tau _r^E\in E^{2r,r}\left( {\rm MGL}\left( r\right) \right) $ is the
universal Thom class of Prop. \ref{thomclasses} and the last map is induced
by the ring structure on ${\Bbb E}$..

\begin{theorem}
\label{thom}Let $\left( {\Bbb E},x_E\right) $ be an oriented ring spectrum, $%
X$ a smooth scheme over $k$ and ${\cal E}\rightarrow X$ be a vector bundle
of rank $r$. Then the Thom map {\rm (\ref{thomiso})} 
\[
\Phi _{{\cal E}}:E^{**}\left( X\right) \longrightarrow E^{*\text{ }+\text{ }%
2r,\,*\text{ }+\text{ }r}\left( {\rm Th}\left( {\cal E}/X\right) \right) 
\]
is an isomorphism.
\end{theorem}

\TeXButton{Proof}{\proof} By \cite{MV}, 3.2.17, we have a canonical cofiber
sequence 
\begin{equation}
{\Bbb P}\left( {\cal E}\right) \longrightarrow {\Bbb P}\left( {\cal E}\oplus 
{\bf 1}\right) \longrightarrow {\rm Th}\left( {\cal E}/X\right) .
\label{pthom}
\end{equation}
The projective bundle structure theorem for $E$-cohomology (\cite{Mo},
3.2.10, (1)) together with the $E$-cohomology long exact sequence associated
to (\ref{pthom}) yield a short exact sequence 
\begin{equation}
0\rightarrow E^{**}\left( {\rm Th}\left( {\cal E}/X\right) \right)
\longrightarrow E^{**}\left( X\right) \left[ t\right] /\left( t\cdot f\left(
t\right) \right) \stackrel{\pi }{\longrightarrow }E^{**}\left( X\right)
\left[ s\right] /\left( f\left( s\right) \right) \rightarrow 0  \label{exact}
\end{equation}
where $f\left( t\right) \doteq t^r+c_1\left( {\cal E}\right) t^{r-1}+\cdots
+c_r\left( {\cal E}\right) $ ($t$ corresponding to the first Chern class of
the canonical line bundle on ${\Bbb P}\left( {\cal E}\right) $), ${\bf 1}$
denotes the trivial line bundle over $X$, $\pi $ maps $t$ to $s$ and $%
c_i\left( {\cal E}\right) \in E^{2i,i}\left( X\right) $ denotes the $i$-th
Chern class of the vector bundle ${\cal E}$. We have used that $c_i\left( 
{\cal E}\oplus 1\right) =c_i\left( {\cal E}\right) $ . Then $f\left(
t\right) $ has bidegree $\left( 2n,n\right) $ and is in the kernel of $\pi $%
. Then, there is a unique Thom class $\tau _r^E\left( {\cal E}\right) $ for $%
{\cal E}$ in $E^{2r,r}\left( {\rm Th}\left( {\cal E}/X\right) \right) $
mapping to $f\left( t\right) $. Obviously we have $\tau _r^E\left( {\cal E}%
\right) =\tau _r^E\circ \Sigma ^\infty \left( \lambda _{{\cal E}}\right) $
and the theorem follows from the exactness of (\ref{exact}). 
\TeXButton{End Proof}{\endproof}

\begin{remark}
{\rm As clear from the above proof, the Thom isomorphism for oriented ring
spectra reduces to the projective bundle theorem together with the fact that
''orientability'' of the spectrum implies the existence of ''universal Thom
classes'' which in its turn implies ''orientability'' of any vector bundle.}

{\rm For the construction of }$BP${\rm -spectra we will only need the Thom
isomorphism for the }${\Bbb MGL}_{(p)}${\rm \ of }${\Bbb MGL}${\rm , at a
prime }$p${\rm \ and this actually follows from the Thom isomorphism for }$%
{\Bbb MGL}${\rm . A proof of this case can be found in \cite{Vo2}, Lecture 3.%
}
\end{remark}

\begin{corollary}
\label{stablethom}If $\left( {\Bbb E},x_E\right) $ is an oriented ring
spectrum, there is a canonical Thom isomorphism 
\[
\Phi :E^{**}\left( {\rm BGL}\right) \widetilde{\longrightarrow }E^{**}\left( 
{\Bbb MGL}\right) .
\]
Moreover $\Phi $ restricts, in bidegree $(0,0)$, to a bijection between ring
spectra maps $\Sigma ^\infty BGL_{+}\rightarrow {\Bbb E}$ and ring spectra
maps ${\Bbb MGL}\rightarrow {\Bbb E}$.
\end{corollary}

\TeXButton{Proof}{\proof} The first assertion is just the stable version of
Theorem \ref{thom} applied to the canonical $n$-plane bundles $\xi
_n\rightarrow {\rm BGL}_n$. In fact, the naturality of the Thom diagonal (%
\ref{diagonal}) implies the commutativity of 
\[
\begin{tabular}{lll}
${\Bbb P}^1\wedge {\rm MGL}\left( n\right) $ & $\stackrel{\Delta _{\xi
_n\oplus 1}}{\longrightarrow }$ & $\left( {\rm BGL}_n\right) _{+}\wedge 
{\Bbb P}^1\wedge {\rm MGL}\left( n\right) $ \\ 
$\qquad ^{\sigma _n}\downarrow $ &  & $\qquad \qquad \downarrow ^{\left(
i_n\right) _{+}\wedge \sigma _n}$ \\ 
{\rm MGL}$\left( n+1\right) $ & $\stackunder{\Delta _{\xi _n}}{%
\longrightarrow }$ & $\left( {\rm BGL}_{n+1}\right) _{+}\wedge {\rm MGL}%
\left( n+1\right) $%
\end{tabular}
\]
for any $n>0$, where the $\sigma _n$'s are the bonding maps of algebraic
cobordism, the $i_n$'s are the natural inclusions ${\rm BGL}%
_n\hookrightarrow {\rm BGL}_{n+1}$ and we used that the multiplicativity
property of Thom spaces 
\[
{\rm Th}\left( \left( {\cal E}\oplus {\bf 1}\right) /X\right) \simeq {\Bbb P}%
^1\wedge {\rm Th}\left( {\cal E}/X\right) 
\]
which holds for any vector bundle ${\cal E}$ over $X$. Therefore, for any $%
\left( p,q\right) $, the diagram 
\[
\begin{tabular}{ccc}
$E^{p,q}\left( {\rm BGL}_{n+1}\right) $ & $\stackrel{\left( i_n\right)
_{+}^{*}}{\longrightarrow }$ & $E^{p,q}\left( {\rm BGL}_n\right) $ \\ 
$^{\Phi _{\xi _{n+1}}}\downarrow $ &  & $\downarrow ^{\Phi _{\xi _n}}$ \\ 
$E^{p+2n+2,q+n+1}\left( {\rm MGL}\left( n+1\right) \right) $ &  & $%
E^{p+2n,q+n}\left( {\rm MGL}\left( n\right) \right) $ \\ 
$\qquad \qquad \qquad \qquad \qquad \quad _{\sigma _n^{*}}\searrow $ &  & $%
\nearrow _{\left( {\Bbb P}^1\right) ^{\wedge -1}}$ \\ 
& $E^{p+2n+2,q+n+1}\left( {\Bbb P}^1\wedge {\rm MGL}\left( n\right) \right) $
& 
\end{tabular}
\]
is commutative and so the family of unstable Thom isomorphism $\left( \Phi
_{\xi _n}\right) _{n>0}$ stabilizes to an isomorphism $\Phi :E^{**}\left( 
{\rm BGL}\right) \widetilde{\longrightarrow }E^{**}\left( {\Bbb MGL}\right) $%
.

The second assertion is a long but straightforward verification using the
commutativity of the diagram 
\[
\begin{tabular}{ccc}
$E^{**}\left( {\Bbb MGL}\right) $ & $\stackrel{\mu ^{*}}{\longrightarrow }$
& $E^{**}\left( {\Bbb MGL}\wedge {\Bbb MGL}\right) $ \\ 
$^\Phi \uparrow $ &  & $\uparrow ^{\Phi ^{^{\prime }}}$ \\ 
$E^{**}\left( {\rm BGL}\right) $ & $\stackunder{m^{*}}{\longrightarrow }$ & $%
E^{**}\left( {\rm BGL}\times {\rm BGL}\right) $%
\end{tabular}
\]
where $\mu :{\Bbb MGL}\wedge {\Bbb MGL}\rightarrow {\Bbb MGL}$ is the
product, $m:{\rm BGL}\times {\rm BGL}\rightarrow {\rm BGL}$ is the canonical
map induced by the map ${\rm BGL}_n\times {\rm BGL}_m\rightarrow {\rm BGL}%
_{n+m}$ (and making ${\rm BGL}$ into an $H$-space) and $\Phi ^{\prime }$ are
stable Thom isomorphisms. \TeXButton{End Proof}{\endproof}

\smallskip 

\section{Orientations, ring spectra maps and formal group laws}

In this section we establishes the basic correspondence, well known in the
topological complex oriented case, between orientations, maps of ring
spectra and formal group laws. We notice that our proof works also in the
topological case and avoids the use of Quillen result that the complex
cobordism ring ${\rm MU}^{*}$\ is isomorphic to the Lazard ring, a result
that in fact is not known in the case of algebraic cobordism.

\begin{lemma}
\label{basic}\label{stablethom}If $\left( {\Bbb E},x_E\right) $ is an
oriented ring spectrum, an element 
\[
\varphi \in E^{0,0}\left( {\rm BGL}\right) ={\rm Hom}_{{\cal SH}\left(
k\right) }\left( \Sigma ^\infty {\rm BGL}_{+},{\Bbb E}\right) 
\]
is a map of ring spectra iff it corresponds via the isomorphism of
Proposition {\rm \ref{lullo}} to a power series $\widehat{\varphi }$ of the
form 
\[
\widehat{\varphi }\left( x_1,\ldots ,x_n,\ldots \right) =\prod_{i=1}^\infty
h\left( x_i\right) 
\]
with $h\left( t\right) $ a degree zero homogeneous power series of the form $%
1+\alpha _1t+\alpha _2t^2+\cdots $ .
\end{lemma}

\TeXButton{Proof}{\proof} A map $\varphi \in {\rm Hom}_{{\cal SH}\left(
k\right) }\left( \Sigma ^\infty {\rm BGL}_{+},{\Bbb E}\right) $ is a ring
map iff the following diagram commutes 
\[
\begin{tabular}{ccc}
$\Sigma ^\infty BGL_{+}\wedge \Sigma ^\infty BGL_{+}$ & $\longrightarrow $ & 
$\Sigma ^\infty BGL_{+}$ \\ 
$^{\varphi \wedge \varphi }\downarrow $ &  & $\downarrow ^\varphi $ \\ 
${\Bbb E}\wedge {\Bbb E}$ & $\longrightarrow $ & ${\Bbb E}$%
\end{tabular}
\]
Arguing as in \cite{Sw}, 16.47 pp. 404-406, we see that if $g$ is the map 
\[
g:\left( {\Bbb P}^\infty \right) ^\infty \times \left( {\Bbb P}^\infty
\right) ^\infty \longrightarrow \left( {\Bbb P}^\infty \right) ^\infty 
\]
\[
\left( u_1,\ldots ,u_n,\ldots ;v_1,\ldots ,v_n,\ldots \right) \longmapsto
\left( u_1,v_{1,}u_2,v_2,\ldots ,u_n,v_n,\ldots \right) , 
\]
the diagram 
\[
\begin{tabular}{ccc}
$\left( {\Bbb P}^\infty \right) ^\infty \times \left( {\Bbb P}^\infty
\right) ^\infty $ & $\stackrel{g}{\longrightarrow }$ & $\left( {\Bbb P}%
^\infty \right) ^\infty $ \\ 
$^{\theta \times \theta }\downarrow $ &  & $\downarrow ^\theta $ \\ 
${\rm BGL}\times {\rm BGL}$ & $\stackunder{m}{\longrightarrow }$ & ${\rm BGL}
$%
\end{tabular}
\]
is (homotopy) commutative, where $\theta $ is the map defined in Proposition 
\ref{lullo}. Therefore, $\varphi $ is a ring map iff its power series $%
\widehat{\varphi }$ satisfies the relation 
\[
\widehat{\varphi }\left( x_1,\ldots ,x_n,\ldots \right) \cdot \widehat{%
\varphi }\left( y_1,\ldots ,y_n,\ldots \right) =\widehat{\varphi }\left(
x_1,y_1,x_2,y_2,\ldots ,x_n,y_n,\ldots \right) 
\]
and we easily conclude by defining $h\left( t\right) \doteq \widehat{\varphi 
}\left( t,0,\ldots ,0,\ldots \right) $. \TeXButton{End Proof}{\endproof}

\smallskip\ 

Let $\left( {\Bbb E},x_E\right) $ be an oriented ring spectrum and $w:{\Bbb P%
}^\infty \times {\Bbb P}^\infty \rightarrow {\Bbb P}^\infty $ be the
canonical map making ${\Bbb P}^\infty $ into an $H$-space. By Proposition 
\ref{morel} (iii), $w^{*}\left( x_E\right) $ defines a power series ${\rm F}%
_{x_E}$ in $E^{**}\left[ \left[ x_1,x_2\right] \right] $. But since $x_E$
has bidegree $\left( 2,1\right) $, ${\rm F}_{x_E}$ is actually an element of
the subring $E^{*}\left[ \left[ x_1,x_2\right] \right] $, where $E^{*}\doteq
\oplus _iE^{2i,i}$. Recall that $E^{*}$ is a commutative ring (\ref{morel}
(i)).

\begin{proposition}
\label{FGL}If $\left( {\Bbb E},x_E\right) $ is an oriented ring spectrum and 
$w:{\Bbb P}^\infty \times {\Bbb P}^\infty \rightarrow {\Bbb P}^\infty $ the
canonical map making ${\Bbb P}^\infty $ into an $H$-space, the power series $%
{\rm F}_{x_E}\doteq $ $w^{*}\left( x_E\right) $ is a formal group law {\rm (%
\cite{H})} on the ''topological'' subring $E^{*}$.
\end{proposition}

\TeXButton{Proof}{\proof}This is a standard consequence of the fact that $w$
defines an $H$-structure on ${\Bbb P}^\infty $ (see for example \cite{Ru},
VII.6.2). \TeXButton{End Proof}{\endproof}

\begin{theorem}
\label{main} Let $\left( {\Bbb E},x_E\right) $ be an oriented ring spectrum
and ${\rm F}_{x_E}$ the formal group law on $E^{*}$ associated to the given
orientation $x_E$. Then the following sets correspond bijectively:

{\rm (i) }orientations on ${\Bbb E}$ ;

{\rm (ii) }maps of ring spectra ${\Bbb MGL}\rightarrow {\Bbb E}$;

{\rm (iii)} pairs $\left( {\rm F},\varepsilon \right) $ where ${\rm F}$ is a
formal group law on $E^{*}$ and $\varepsilon :{\rm F}\stackrel{\sim }{%
\rightarrow }{\rm F}_{x_E}$ is an isomorphism of formal group laws.
\end{theorem}

\TeXButton{Proof}{\proof} By Proposition \ref{morel} (ii) an orientation $x$
on ${\Bbb E}$ is of the form $f\left( x_E\right) $ where $f\left( t\right)
=t+\alpha _2t^2+\alpha _3t^3+\cdots $ . Such an $f$ gives an isomorphism of
formal group laws ${\rm F}_x\stackrel{\sim }{\rightarrow }{\rm F}_{x_E}$ and
viceversa. Hence (i) and (iii) are in bijection.

With the same notations, the power series 
\[
\widehat{\varphi }\left( x_1,\ldots ,x_n,\ldots \right) \doteq
\prod_{i=1}^\infty \frac{f\left( x_i\right) }{x_i} 
\]
defines a map of ring spectra $\varphi :{\Bbb MGL}$ $\rightarrow {\Bbb E}$,
by Lemma \ref{basic} and Theorem \ref{stablethom}, and this construction can
be inverted. \TeXButton{End Proof}{\endproof}

\begin{remark}
{\rm Replacing }${\cal SH}\left( k\right) ${\rm \ with the topological
stable homotopy category and }${\Bbb MGL}${\rm \ with the complex cobordism
spectrum }${\Bbb MU}${\rm , the proofs of Lemma \ref{basic} and Theorem \ref
{main} carry over without modifications. This is a slightly different
approach with respect to the usual one in topology where one uses Quillen
Theorem (i.e. the isomorphism of MU}$^{*}${\rm \ with the universal Lazard
ring) to prove the equivalence between (ii) and (iii) in Theorem \ref{main}.
Actually, our proof is forced to be different from that since neither }${\rm %
MGL}^{**}${\rm \ nor }${\rm MGL}^{*}${\rm \ are known. Moreover, note that
Quillen Theorem is actually stronger than Theorem \ref{main}, in the
topological case, in the sense that Quillen's result cannot be deduced from
(the topological version of) Theorem \ref{main}.}
\end{remark}

\section{The motivic Quillen idempotent and Brown-Peterson spectra in ${\cal %
SH}\left( k\right) $}

We are going to apply Theorem \ref{main} to the localization of the
algebraic cobordism spectrum at a prime $p$.

Throughout this section, we fix a prime $p$. Let ${\Bbb MGL}_{(p)}$ be the
(Bousfield) localization of ${\Bbb MGL}$ at the prime $p$. Since the
localization map 
\[
\ell :{\Bbb MGL}\longrightarrow {\Bbb MGL}_{(p)} 
\]
is a map of ring spectra, it maps the canonical orientation $x_{{\Bbb MGL}}$
of ${\Bbb MGL}$ (Remark \ref{canonical}) to an orientation $x_{(p)}$ of $%
{\Bbb MGL}_{(p)}$. Hence, ${\Bbb MGL}_{(p)}$ is canonically oriented by $%
x_{(p)}$. Let us denote by {\rm F}$_{x_{(p)}}$ the corresponding formal
group law on ${\rm MGL}_{(p)}^{*}\doteq \oplus _i{\rm MGL}_{(p)}^{2i,i}$
(Proposition \ref{FGL}). Since ${\rm MGL}_{(p)}^{*}$ is a commutative
(Proposition \ref{morel}) ${\Bbb Z}_{(p)}$-algebra, by Cartier theorem,
there exists a canonical strict isomorphism of formal group laws on ${\rm MGL%
}_{(p)}^{*}$ 
\[
\varepsilon :{\rm F}_{x_{(p)}^0}\widetilde{\longrightarrow }{\rm F}%
_{x_{(p)}} 
\]
with ${\rm F}_{x_{(p)}^0}$ a $p$-typical formal group law (\cite{H},
16.4.14). Therefore, by Theorem \ref{main} applied to the oriented ring
spectrum $\left( {\Bbb MGL}_{(p)},x_{(p)}\right) $, to such an $\varepsilon $
is uniquely associated a ring spectra map 
\begin{equation}
{\rm e}:{\Bbb MGL\longrightarrow MGL}_{(p)}.  \label{e}
\end{equation}
If ${\rm C}$ denotes the cofiber of the localization map $\ell :{\Bbb MGL}%
\longrightarrow {\Bbb MGL}_{(p)}$, clearly one has 
\begin{equation}
{\rm MGL}_{(p)}^{**}\left( {\rm C}\right) =0  \label{cofiber}
\end{equation}
and therefore the natural map 
\[
\ell ^{*}:{\rm MGL}_{(p)}^{**}\left( {\Bbb MGL}_{(p)}\right) \longrightarrow 
{\rm MGL}_{(p)}^{**}\left( {\Bbb MGL}\right) 
\]
is an isomorphism.

\begin{proposition}
\label{facile}The isomorphism $\ell ^{*}$ establishes, in bidegree $\left(
0,0\right) $, a bijection between ring spectra maps ${\Bbb MGL}%
_{(p)}\rightarrow {\Bbb MGL}_{(p)}$ and ring spectra maps ${\Bbb MGL}%
\rightarrow {\Bbb MGL}_{(p)}$.
\end{proposition}

\TeXButton{Proof}{\proof}One direction is clear since $\ell $ is a map of
ring spectra. On the other hand, let us consider a map of ring spectra $%
\alpha :{\Bbb MGL}\rightarrow {\Bbb MGL}_{(p)}$ and let $\beta :{\Bbb MGL}%
_{(p)}\rightarrow {\Bbb MGL}_{(p)}$ the unique map such that $\ell
^{*}\left( \beta \right) =\alpha $. We must prove that $\beta $ is a map of
ring spectra.

In the following diagram 
\[
\begin{tabular}{ccc}
${\Bbb MGL}\wedge {\Bbb MGL}$ & $\stackrel{\mu }{\longrightarrow }$ & ${\Bbb %
MGL}$ \\ 
$^{\ell \wedge \ell }\downarrow $ &  & $\downarrow ^\ell $ \\ 
${\Bbb MGL}_{(p)}\wedge {\Bbb MGL}_{(p)}$ & $\stackrel{\mu _{(p)}}{%
\longrightarrow }$ & ${\Bbb MGL}_{(p)}$ \\ 
$^{\beta \wedge \beta }\downarrow $ &  & $\downarrow ^\beta $ \\ 
${\Bbb MGL}_{(p)}\wedge {\Bbb MGL}_{(p)}$ & $\stackunder{\mu _{(p)}}{%
\longrightarrow }$ & ${\Bbb MGL}_{(p)}$%
\end{tabular}
\]
(where the horizontal arrows are product maps) the upper square is
commutative and the outer square too since $\alpha $ is a map of ring
spectra. If $d$ denotes the difference 
\[
\beta \circ \mu _{(p)}-\mu _{(p)}\circ \left( \beta \wedge \beta \right) 
\text{ ,} 
\]
we know that $d\circ \left( \ell \wedge \ell \right) $ is zero. But since 
\[
\left( {\Bbb MGL}\wedge {\Bbb MGL}\right) _{(p)}\simeq {\Bbb MGL}%
_{(p)}\wedge {\Bbb MGL}_{(p)} 
\]
and 
\[
{\rm MGL}_{(p)}^{**}\left( {\Bbb MGL\wedge MGL}\right) \longrightarrow {\rm %
MGL}_{(p)}^{**}(\left( {\Bbb MGL}\wedge {\Bbb MGL}\right) _{(p)}) 
\]
is an isomorphism by the same argument used in (\ref{cofiber}), we conclude
that also $d$ is zero i.e. that $\beta $ is indeed a map of ring spectra. 
\TeXButton{End Proof}{\endproof}

\begin{corollary}
The unique map {\rm e}$_{(p)}:{\Bbb MGL}_{(p)}\rightarrow {\Bbb MGL}_{(p)}$
such that $\ell ^{*}({\rm e}_{(p)})=$ {\rm e}, is a map of ring spectra.
\end{corollary}

Since the canonical procedure to make a given formal group law $p$-typical,
is trivial when applied to a formal group law which is already $p$-typical (%
\cite{H}, 31.1.9, p. 429), the ring map {\rm e}$_{(p)}$ is idempotent. We
call {\rm e}$_{(p)}$ the {\em motivic Quillen idempotent}.

\begin{definition}
The {\rm Brown-Peterson spectrum} in ${\cal SH}\left( k\right) $ associated
to the prime $p$ is the spectrum ${\Bbb BP}$ colimit of the diagram of ring
spectra and ring spectra maps in ${\cal SH}\left( k\right) $%
\[
\cdots \rightarrow {\Bbb MGL}_{(p)}\stackrel{{\rm e}_{(p)}}{\rightarrow }%
{\Bbb MGL}_{(p)}\stackrel{{\rm e}_{(p)}}{\rightarrow }{\Bbb MGL}_{(p)}%
\stackrel{{\rm e}_{(p)}}{\rightarrow }{\Bbb MGL}_{(p)}\rightarrow \cdots 
\]
\end{definition}

Therefore, ${\Bbb BP}$ is a commutative ring spectrum and there are
canonical maps of ring spectra ${\rm u}:{\Bbb BP\rightarrow MGL}_{(p)}$ and $%
\widetilde{{\rm e}}:{\Bbb MGL}_{(p)}\rightarrow {\Bbb BP}$ such that $%
\widetilde{{\rm e}}\circ {\rm u}={\rm id}_{{\Bbb BP}}$ and ${\rm u}\circ 
\widetilde{{\rm e}}={\rm e}_{(p)}$. In particular, ${\Bbb BP}$ is a direct
summand of ${\Bbb MGL}_{(p)}$.

\begin{remark}
{\rm Note that we were forced (unlike in the topological case) to prove the
existence of the Quillen idempotent without resorting to Quillen theorem
which is not known to hold in the algebraic case. However, the construction
of the idempotent given above works in the topological case too, hence
yielding a different construction from the usual one that uses the
isomorphism between }${\rm MU}^{*}${\rm \ and the Lazard ring.}
\end{remark}

\medskip\

\smallskip\ 

E-mail: {\tt vezzosi@dm.unibo.it}


\begin{thebibliography}{99}
\bibitem{CGK}  M. Cole, J.P.C. Greenless, I. Kriz, {\em The universality of
equivariant complex bordism}, preprint, 1999,
http://hopf.math.purdue.edu/cgi-bin/generate?/pub/Cole-Greenlees-Kriz/AThom

\bibitem{Fu}  W. Fulton, {\em Intersection theory}, Second edition, Springer
Verlag, New York, 1998.

\bibitem{H}  M. Hazewinkel, {\em Formal groups and applications}, Academic
Press, New York, 1978.

\bibitem{Mo}  F. Morel, {\em Basic properties of the stable homotopy
categories of smooth schemes}, Preprint,
http://www.math.jussieu.fr/~morel/Satble.ps

\bibitem{MV}  F. Morel, V. Voevodsky, ${\Bbb A}^1${\em -homotopy of schemes}%
, to appear in Publ. Math. IHES.

\bibitem{Q1}  D. Quillen, {\em On the formal group laws of unoriented and
complex cobordism theory}, Bull. A.M.S. {\bf 75}, 1969, 1293-1298.

\bibitem{Q2}  D. Quillen, {\em Elementary proofs of some results of
cobordism theory using Steenrod operations}, Adv. in Math., 7, 1971, 29-56.

\bibitem{Ru}  Y. B. Rudyak, {\em On Thom spectra, orientability and cobordism%
}, Springer Monographs in Mathematics, Springer Verlag, New York, 1998.

\bibitem{St}  R. Stong, {\em Notes on cobordism theory}, Princeton
University Press, Princeton, 1968.

\bibitem{Sw}  R. M. Switzer, {\em Algebraic topology-Homotopy and Homology},
Springer Verlag, New York, 1975.

\bibitem{To}  B. Totaro, {\em Torsion algebraic cycles and complex cobordism}%
, Journal A.M.S., {\bf 10}, 1992, 467-493.

\bibitem{Vo1}  V. Voevodsky, ${\Bbb A}^1${\em -homotopy theory}, Doc. Math.
J. DMV, Extra volume ICM 1998, I, 417-442.

\bibitem{Vo2}  V. Voevodsky, {\em Lecture 3: Algebraic cobordism}, Workshop
on homotopy theory of algebraic varieties, MSRI, Berkeley, May 1998,
http://www.msri.org/publications/ln/msri/ 1998/homotopy/voevodsky/3.

\bibitem{Vo3}  V. Voevodsky, {\em The Milnor conjecture}, preprint MPI, 1996.

\bibitem{We}  C. Weibel, {\em Voevodsky's Seattle Lectures: }$K${\em -theory
and motivic cohomology}, Proc. symp. Pure Math. {\bf 67}, 1999. 

\bibitem{W}  W. S. Wilson, {\em The }$\Omega ${\em -spetrum for
Brown-Peterson cohomology, part II}, Am. J. Math., 97 (1975), 101-123.
\end{thebibliography}
\end{document}